\documentclass[12pt]{article}
\usepackage{latexsym,amssymb,amsmath}
\begin{document}
\baselineskip=18pt

\newcommand{\la}{\langle}
\newcommand{\ra}{\rangle}
\newcommand{\psp}{\vspace{0.4cm}}
\newcommand{\pse}{\vspace{0.2cm}}
\newcommand{\ptl}{\partial}
\newcommand{\dlt}{\delta}
\newcommand{\sgm}{\sigma}
\newcommand{\al}{\alpha}
\newcommand{\be}{\beta}
\newcommand{\G}{\Gamma}
\newcommand{\gm}{\gamma}
\newcommand{\vs}{\varsigma}
\newcommand{\Lmd}{\Lambda}
\newcommand{\lmd}{\lambda}
\newcommand{\td}{\tilde}
\newcommand{\vf}{\varphi}
\newcommand{\yt}{Y^{\nu}}
\newcommand{\wt}{\mbox{wt}\:}
\newcommand{\rd}{\mbox{Res}}
\newcommand{\ad}{\mbox{ad}}
\newcommand{\stl}{\stackrel}
\newcommand{\ol}{\overline}
\newcommand{\ul}{\underline}
\newcommand{\es}{\epsilon}
\newcommand{\dmd}{\diamond}
\newcommand{\clt}{\clubsuit}
\newcommand{\vt}{\vartheta}
\newcommand{\ves}{\varepsilon}
\newcommand{\dg}{\dagger}
\newcommand{\tr}{\mbox{Tr}}
\newcommand{\ga}{{\cal G}({\cal A})}
\newcommand{\hga}{\hat{\cal G}({\cal A})}
\newcommand{\Edo}{\mbox{End}\:}
\newcommand{\for}{\mbox{for}}
\newcommand{\kn}{\mbox{ker}}
\newcommand{\Dlt}{\Delta}
\newcommand{\rad}{\mbox{Rad}}
\newcommand{\rta}{\rightarrow}
\newcommand{\mbb}{\mathbb}

\begin{center}{\LARGE \bf Poisson and Hamiltonian Superpairs  over}\end{center}
\begin{center}{\LARGE \bf   Polarized Associative Algebras}
\end{center}
\vspace{0.2cm}

\begin{center}{\large Xiaoping Xu}\end{center}
\begin{center}{Department of Mathematics, The Hong Kong University of Science \& Technology}\end{center}
\begin{center}{Clear Water Bay, Kowloon, Hong Kong}\footnote{Research supported
 by Hong Kong RGC Competitive Earmarked Research Grant HKUST6133/00P.}\end{center}

\vspace{0.3cm}

\begin{center}{\Large \bf Abstract}\end{center}
\vspace{0.2cm}

{\small Poisson superpair is a pair of Poisson superalgebra structures on a super commutative associative algebra, whose any linear combination is also a Poisson superalgebra structure. In this paper, we first construct certain linear and quadratic Poisson superpairs over a  semi-finitely-filtered polarized $\Bbb{Z}_2$-graded associative algebra. Then we give a construction of certain Hamiltonian superpairs in the formal variational calculus over any finite-dimensional $\Bbb{Z}_2$-graded associative algebra with a supersymmetric nondegenerate associative bilinear form. Our constructions are based on the Adler mapping in a general sense. Our works in this paper can be viewed as noncommutative generalizations of the Adler-Gel'fand-Dikii Hamiltonian pair}.  
\vspace{0.2cm}

\section{Introduction}

In the theory of completely integrable systems, one of the most beautiful structures is the Adler-Gel'fand-Dikii Hamiltonian pair, which was constructed through the Adler mapping (cf. [A], [GDi3]). Such a Hamiltonian pair gives a pair of Poisson structures on the quotient space of the differential polynomial algebra of scalar differential operators with fixed order modulo its subspace of total differential polynomials (flux). They do not form Poisson algebra structures because the quotient space does not form an associative algebra. 

Gel'fand and Dorfman [GDo2] generalized the Adler-Gel'fand-Dikii Hamiltonian pair to that over differential operators with fixed order and the coefficients in a matrix algebra. The first Hamiltonian structure for the Kadomtsev-Petviashvili hierarchy was suggested by Watanabe [W]. Dickey [D1] added the second Hamiltonian structure for the Kadomtsev-Petviashvili hierarchy. These two structures form an analogue of the Adler-Gel'fand-Dikii Hamiltonian pair over scalar pseudo-differential operators of positive order one and  infinite negative order. Radul [R] generalized it over scalar pseudo-differential operators of finite positive order and infinite negative order.

Manin and Radul [MR] gave a supersymmetric extension of the Kadomtsev-Petviashvili hierarchy.  Das and Huang [DH] essentially partially generalized the Adler-Gel'fand-Dikii's construction over  scalar differential operators with fixed order to that over scalar pseudo-differential operators with fixed positive and negative orders. Dofman and Fokas [DF] generalized
the Adler-Gel'fand-Dikii Hamiltonian pair to that over differential operators with fixed order and the coefficients in an algebra of pseudo-differential operators, in which the algebra plays the same role as a number field. It seems that there is a problem of how to interpret Dorfman and Fokas' results in [DF] in terms of the theory of Hamiltonian operators over a field. 

In [X1], we generalized the theory of Hamiltonian operators to that of Hamiltonian superoperators over fermionic fields. Moreover, we established  an analogous theory of supervariables in [X2]. In [X3], we proved that conformal superalgebras are equivalent to certain linear Hamiltonian superoperator of super functions in one real variable. 

We observe that the Adler-Gel'fand-Dikii's construction was essentially based on a polarization of the algebra of pseudo-differential operators. Their construction could be generalized and applied to
more general polarized associative algebras. The first objective in this paper is to construct linear and quadratic Poisson superpairs on supersymmetric polynomial functions of  a  semi-finitely-filtered polarized $\Bbb{Z}_2$-graded associative algebra. The second objective is to construct certain Hamiltonian superpairs in the formal variational calculus over any finite-dimensional $\Bbb{Z}_2$-graded associative algebra with a supersymmetric nondegenerate associative bilinear form. Our constructions are based on the Adler mapping in a general sense. The results in this paper tell that there is a deep algebraic essence behind the  Adler-Gel'fand-Dikii Hamiltonian pair. In fact, the linear structure in our general  analugues of the  Adler-Gel'fand-Dikii Hamiltonian pair depneds on a central element. 

We shall give a more technical introduction in Section 2. In Section 3, we shall present certain structural properties and constructions of $\Bbb{Z}_2$-graded polarized associative algebras. The  Poisson superpairs will be given in Section 4. Section 5 is devoted to the Hamiltonian superpairs in the formal variational calculus.

\section{Technical Background}

This section serves as a technical introduction of the whole paper. 

Throughout this paper, we let $\Bbb{F}$ a field with characteristic not equal to two unless it is specified. All the vector spaces (algebras) are assumed over $\Bbb{F}$. Denote by $\Bbb{Z}$ the ring of integers and by $\Bbb{N}$ the set of nonnegative integers.
For any two integers $m_1,m_2$,  we shall often use the following notation of index throughout this paper:
$$\ol{m_1,m_2}=\left\{\begin{array}{ll}\{m_1,m_1+1,m_1+2,...,m_2\}&\mbox{if}\;\;m_1\leq m_2,\\\emptyset&\mbox{if}\;\;m_1>m_2.\end{array}\right.\eqno(2.1)$$

First we introduce the definition of abstract Hamiltonian superoperators.  Let $({\cal G}, [\cdot,\cdot])$ be a Lie superalgebra and let $M$ be a ${\cal G}$-module. 
For a positive integer $q$,  a $q$-{\it form of} ${\cal G}$ {\it with values in} $M$ is a multi-linear map $\omega:\;{\cal G}^q={\cal G}\times \cdots \times {\cal G}\rightarrow M$ for which
$$\omega (\xi_1,\xi_2, \cdots,\xi_q)=-(-1)^{ij}\omega(\xi_1,\cdots,\xi_{\ell-1},\xi_{\ell+1},\xi_{\ell},\xi_{\ell+2},\cdots, \xi_q)\eqno(2.2)$$
for $\xi_k\in{\cal G}, \xi_{\ell}\in {\cal G}_i$ and $\xi_{\ell+1}\in {\cal G}_j.$
We denote by $c^q({\cal G},M)$ the set of $q$-forms. Moreover, we define a differential $d:\;c^q({\cal G},M)\rightarrow c^{q+1}({\cal G},M)$ by
\begin{eqnarray*} d\omega(\xi_1,\xi_2,...,\xi_{q+1})&=&\sum_{\ell=1}^{q+1}(-1)^{\ell+1+i_{\ell}(i_1+\cdots i_{\ell-1})}\xi_{\ell}\omega(\xi_1,...,\check{\xi}_{\ell},...,\xi_{q+1})\\& &+\sum_{\ell_1<\ell_2}(-1)^{\ell_1+\ell_2+(i_{\ell_1}+i_{\ell_2})(i_1+\cdots+i_{\ell_1-1})+i_{\ell_2}(i_{\ell_1+1}+\cdots+i_{\ell_2-1})}\\& &\omega([\xi_{\ell_1},\xi_{\ell_2}],\xi_1,...,\check{\xi}_{\ell_1},...,\check{\xi}_{\ell_2},....,\xi_{q+1})\hspace{4.6cm}(2.3)\end{eqnarray*}
for $\omega\in c^q({\cal G},M)$ and $\xi_k\in {\cal G}_{i_k}$ with $k\in\ol{1,q+1}$, where the above index ``check'' means deleting the term under it. A $q$-form $\omega$ is called {\it closed} if $d\omega=0$.

For any $u\in M$, we define a one-form $du$ by
$$du(\xi)=\xi(u)\qquad\for\;\;\xi\in{\cal G}.\eqno(2.4)$$
Let $\Omega$ be a subspace of $c^1({\cal G},M)$ such that $dM\subset \Omega$. Suppose that $H:\;\Omega \rightarrow {\cal G}$ is a linear map. We call $H$ $\Bbb{Z}_2$-{\it graded} if
$$H(\Omega)=H(\Omega)_0\oplus H(\Omega)_1,\qquad\mbox{where}\qquad H(\Omega)_i=H(\Omega)\bigcap{\cal G}_i.\eqno(2.5)$$
Moreover, $H$ is called {\it super skew-symmetric} if
$$\phi_1(H\phi_2)=-(-1)^{i_1i_2}\phi_2(H\phi_1)\qquad \mbox{where}\qquad H\phi_j\in H(\Omega)_{i_j}.\eqno(2.6)$$
For a $\Bbb{Z}_2$-graded super skew-symmetric linear map $H: \Omega\rightarrow {\cal G}$, we define a 2-form $\omega_H$ on $H(\Omega)$ by 
$$\omega_H(H\phi_1,H\phi_2)=\phi_2(H\phi_1)\qquad\for\;\;\phi_1,\phi_2\in \Omega. \eqno(2.7)$$

 We say that a super skew-symmetric $\Bbb{Z}_2$-graded linear map $H: \Omega\rightarrow {\cal G}$ is a {\it Hamiltonian superoperator}  if

(a) the subspace $H(\Omega)$ of ${\cal G}$  forms a subalgebra;

(b) the left and right radicals of the form $\omega_H$  are $\Bbb{Z}_2$-graded and $d\omega_H\equiv 0$  on $H(\Omega)$.

Two $\Bbb{Z}_2$-graded linear maps $H_1, H_2: \Omega\rightarrow {\cal G}$ are called a {\it Hamiltonian pair} if $\lmd_1H_1+\lmd_2H_2$ is a Hamiltonian superoperator for any $\lmd_1,\lmd_2\in\Bbb{F}$.

Next we introduce the Adler-Gel'fand-Dikii Hamiltonian pair and the known generalizations. We assume that  $\Bbb{F}$ is a field of real numbers or a field of complex numbers. Let $k$ be a positive integer and let $\{u_0,u_1,...,u_{k-1}\}$ be $k$ $C^{\infty}$-functions in the real variable $x$. Set
$$u_j^{(m)}={d^mu_j\over dx^m}\qquad\for\;\;m\in\Bbb{N},\;j\in\ol{0,k-1}.\eqno(2.8)$$
Denote
$${\cal P}=\Bbb{F}[u_j^{(m)}\mid m\in\Bbb{N},\;j\in\ol{0,k-1}\}],\eqno(2.9)$$
the differential polynomial algebras of $\{u_j(x)\mid j\in\ol{0,k-1}\}$. We view
$${d\over dx}=\sum_{j\in\ol{0,k-1},\;m\in\Bbb{N}}u_j^{(m+1)}\ptl_{u_j^{(m)}}\eqno(2.10)$$
as a derivation of ${\cal P}$. For convenience, we denote
$$\ptl={d\over dx}\eqno(2.11)$$
and define the algebra of pseudo-differential operators
$${\cal D}=\{\sum_{l=-\infty}^nf_l\ptl^l\mid n\in\Bbb{Z},\;f_l\in{\cal P}\}\eqno(2.12)$$
with the multiplication determined by
$$(f\ptl^m)(g\ptl^n)=\sum_{p=0}^{\infty}\left(\!\!\begin{array}{c}m \\ p\end{array}\!\!\right)fg^{(p)}\ptl^{m+n-p}\qquad\for\;\;f,g\in{\cal P},\;\;m,n\in\Bbb{N},\eqno(2.13)$$
where 
$$g^{(p)}=\ptl^p(g).\eqno(2.14)$$
 Set 
$${\cal G}=\sum_{j=0}^{k-1}{\cal P}\ptl^j\subset{\cal D}\eqno(2.15)$$
and
$$\ptl_X=\sum_{j\in\ol{0,k-1},\;m\in\Bbb{N}}a_j^{(m)}\ptl_{u^{(m)}_j}\qquad\for\;\;X=\sum_{j=0}^{k-1}a_j\ptl^j\in {\cal G}.\eqno(2.16)$$
As derivations of ${\cal P}$, 
$$[\ptl_X,\ptl]=0\qquad\for\;\;X\in{\cal G}.\eqno(2.17)$$
Define
$$\ptl_X(Y)=\sum_{j=0}^{k-1}\ptl_X(b_j)\ptl^j\qquad\for\;\;X,\;Y=\sum_{j=0}^{k-1}b_j\ptl^j\in {\cal G}.\eqno(2.18)$$
The Lie bracket on ${\cal G}$ is defined by
$$[X,Y]_0=\ptl_X(Y)-\ptl_Y(X)\qquad\for\;\;X,Y\in{\cal G}.\eqno(2.19)$$

Set
$$\td{\cal P}={\cal P}/\ptl({\cal P})\eqno(2.20)$$
and use notation
$$\td{f}=f+\ptl({\cal P})\qquad\for\;\;f\in{\cal P}.\eqno(2.21)$$
Moreover, by (2.17),  we define an action of ${\cal G}$ on $\td{\cal P}$:
$$X(\td{f})=(\ptl_X(f))^{\sim}\qquad\for\;\;X\in{\cal G},\;f\in{\cal A}.\eqno(2.22)$$
 Then $\td{\cal P}$ forms a ${\cal G}$-module. 

Define
$$\Omega=\sum_{j=0}^{k-1}\ptl^{-1-j}{\cal P}\subset{\cal D}\eqno(2.23)$$
and identify it with a subspace of one-forms by
$$\xi(X)=\sum_{j=0}^{k-1}(a_i\al_i)^{\sim}\qquad\for\;\;\xi=\sum_{j=0}^{k-1}\ptl^{-1-j}\al_j\in\Omega,\;X=\sum_{j=0}^{k-1}a_j\ptl^j\in {\cal G}.\eqno(2.24)$$
Moreover, we define the projection from ${\cal D}$ to ${\cal G}$ by
$$(\sum_{l=-\infty}^nf_l\ptl^l)_+=\sum_{l=0}^nf_l\ptl^l\qquad\for\;\;\sum_{l=-\infty}^nf_l\ptl^l\in {\cal D}\;\;\mbox{with}\;\;n\in\Bbb{N}.\eqno(2.25)$$
Set
$$L=\ptl^k+\sum_{j=0}^{k-1}u_j\ptl^j.\eqno(2.26)$$
We define two linear maps $H_1,H_2: \Omega \rta {\cal G}$ by
$$H_1(\xi)=([L,\xi])_+,\;\;H_2(\xi)=(L\xi)_+L-L(\xi L)_+.\eqno(2.27)$$
Then $H_1$ and $H_2$ forms a Hamiltonian pair, which is called the {\it  Adler-Gel'fand-Dikii Hamiltonian pair} (cf. [A], [GD3], [D]). The map $H_2$ is called the {\it Adler mapping} (cf. [A]).
 
Gel'fand and Dorfman [GDo2] generalized the Adler-Gel'fand-Dikii Hamiltonian pair over the differential operator $L$ with $u_j$ taking values in an $m\times m$ matrix algebra.  Watanabe [W] and Dickey [D1] obtained Adler-Gel'fand-Dikii Hamiltonian pairs over the pseudo-differential operator
$$L=\ptl+\sum_{j=-\infty}^0u_j\ptl^j.\eqno(2.28)$$
Radul [R] further generalized it over the pseudo-differential operator
$$L=\ptl^k+\sum_{j=-\infty}^ku_j\ptl^j\eqno(2.29)$$
for any positive integer $k$. Das and Huang [DH] generalized the Gel'fand-Dikii's construction of Hamiltonian pairs essentially by extending ${\cal P}$ to the algebra of differential polynomials of $2k$ $C^{\infty}$-functions $\{u_{-k},...,u_{-1}, u_0,...,u_{k-1}\}$, taking
$$L=\ptl^k+\sum_{j=-k}^{k-1}u_j\ptl^j\eqno(2.30)$$
and keeping ${\cal G}$ and $\Omega$ the same as in (2.15) and (2.23) with the new ${\cal P}$. Dorfman and Fokas obtained the Adler-Gel'fand-Dikii Hamiltonian pair over the differential operator $L$ with $u_j$ taking values in an algebra of pseudo-differential operators, in which the algebra plays the same role as a number field. It seems that there is a problem of how to interpret Dorfman and Fokas' results in [DF] in terms of the theory of Hamiltonian operators over a field. 

Let $\Bbb{F}$ be a general field. A $\Bbb{Z}_2$-graded associative algebra 
$${\cal A}={\cal A}_0\oplus{\cal A}_1\eqno(2.31)$$
is called a {\it polarized} $\Bbb{Z}_2$-{\it graded associative algebra} if
${\cal A}$ has a nondegenerate $\Bbb{Z}_2$-graded supersymmetric associative  bilinear form $\la\cdot,\cdot\ra:{\cal A}\times {\cal A}\rta \Bbb{F}$, that is,
$$\la{\cal A}_0,{\cal A}_1\ra=\{0\},\;\;\la u,v\ra=(-1)^{i_1i_2}\la v,u\ra,\;\;\la uv,w\ra=\la u,vw\ra\eqno(2.32)$$
for $u\in{\cal A}_{i_1},\;v\in{\cal A}_{i_2},\;w\in{\cal A}$, and  contains two $\Bbb{Z}_2$-graded isotropic subalgebras ${\cal A}^+$ and ${\cal A}^-$, namely,
$$\la {\cal A}^+,{\cal A}^+\ra=\{0\},\;\;\la {\cal A}^-,{\cal A}^-\ra=\{0\},\eqno(2.33)$$
 such that 
$${\cal A}={\cal A}^+\oplus {\cal A}^-.\eqno(2.34)$$
Expressions (2.33) and (2.34) are called a {\it polarization of} ${\cal A}$. 

An associative algebra ${\cal B}$ is called {\it super commutative} if ${\cal B}={\cal B}_0\oplus {\cal B}_1$ is a $\Bbb{Z}_2$-graded algebra such that
$$uv=(-1)^{i_1i_2}vu\qquad\for\;\;u\in {\cal B}_{i_1},\;v\in{\cal B}_{i_2}.\eqno(2.35)$$
A {\it Poisson superalgebra} is a super commutative associative ${\cal B}$ with another algebraic operation $\{\cdot,\cdot\}$, called {\it Poisson superbracket}, such that $({\cal B},\{\cdot,\cdot\})$ forms a Lie superalgebra and the following compatibility condition holds
$$\{u,vw\}=\{u,v\}w+(-1)^{i_2i_2}v\{u,w\}\qquad\for\;\;u\in{\cal B}_{i_1},\;v\in{\cal B}_{i_2},\;w\in{\cal B}.\eqno(2.36)$$

The main purpose of this paper is to give a more systematic study in generalizations of the  Adler-Gel'fand-Dikii's construction of  Hamiltonian pairs over a  polarized $\Bbb{Z}_2$-graded associative algebra.
In particular, we obtain  pairs of a linear and  a quadratic Poisson superalgebra structures on supersymmetric polynomial functions of  a  semi-finitely-filtered polarized $\Bbb{Z}_2$-graded associative algebras. 

\section{Polarized Associative Algebras}

In this section, we shall present certain structural properties and constructions of $\Bbb{Z}_2$-graded polarized associative algebras.

As we shall show below, the structure of a $\Bbb{Z}_2$-graded polarized associative algebra is determined by a certain compatible pair of $\Bbb{Z}_2$-graded associative algebra structures on a vector space. 

 Let 
$${\cal A}={\cal A}^+\oplus {\cal A}^-\eqno(3.1)$$
be a polarized $\Bbb{Z}_2$-graded associative algebra with the nondegenerate $\Bbb{Z}_2$-graded supersymmetric associative  bilinear form $\la\cdot,\cdot\ra$. Set
$${\cal A}^{\pm}_i={\cal A}_i\bigcap {\cal A}^{\pm}\qquad\for\;\;i\in\Bbb{Z}_2.\eqno(3.2)$$
Take a basis $\{\vs_{i,j}^+\mid j\in I_i\}$ of ${\cal A}^+_i$ for $i\in\Bbb{Z}_2$, where $I_i$ are index sets. Suppose that ${\cal A}^-_i$ has a dual basis $\{\vs_{i,j}^-\mid j\in I_i\}$ with respect to the basis $\{\vs_{i,j}^+\mid j\in I_i\}$ of ${\cal A}^+_i$ for $i\in\Bbb{Z}_2$, that is,
$$\la \vs_{i_1,j_1}^-,\vs_{i_2,j_2}^+\ra=\dlt_{i_1,i_2}\dlt_{j_1,j_2}\qquad\for\;\;i_1,i_2\in\Bbb{Z}_2,\;j_1\in I_{i_1},\;j_2\in I_{i_2}.\eqno(3.3)$$
This assumption trivially holds when ${\cal A}$ is finite-dimensional.  

Write
$$\vs_{i_1,j_1}^{\pm}\vs_{i_2,j_2}^{\pm}=\sum_{j_3\in I_{i_1+i_2}}a^{\pm,j_3}_{i_1,j_1;i_2,j_2}\vs^{\pm}_{i_1+i_2,j_3}\qquad\for\;\;i_1,i_2\in\Bbb{Z}_2,\;j_1\in I_{i_1},\;j_2\in I_{i_2}.\eqno(3.4)$$
Then
$$\la \vs^+_{i_1,j_1}\vs^-_{i_2,j_2},\vs^-_{i_1+i_2,j_3}\ra=\la \vs^+_{i_1,j_1},\vs^-_{i_2,j_2}\vs^-_{i_1+i_2,j_3}\ra=(-1)^{i_1}a^{-,j_1}_{i_2,j_2;i_1+i_2,j_3},\eqno(3.5)$$
$$\la\vs^+_{i_1+i_2,j_3}, \vs^+_{i_1,j_1}\vs^-_{i_2,j_2}\ra=\la\vs^+_{i_1+i_2,j_3}\vs^+_{i_1,j_1},\vs^-_{i_2,j_2}\ra=(-1)^{i_2}a^{+,j_2}_{i_1+i_2,j_3;i_1,j_1},\eqno(3.6)$$
$$\la \vs^-_{i_1,j_1}\vs^+_{i_2,j_2},\vs^+_{i_1+i_2,j_3}\ra=\la \vs^-_{i_1,j_1},\vs^+_{i_2,j_2}\vs^+_{i_1+i_2,j_3}\ra=a^{+,j_1}_{i_2,j_2;i_1+i_2,j_3},\eqno(3.7)$$
$$\la\vs^-_{i_1+i_2,j_3}, \vs^-_{i_1,j_1}\vs^+_{i_2,j_2}\ra=\la\vs^-_{i_1+i_2,j_3}\vs^-_{i_1,j_1},\vs^+_{i_2,j_2}\ra=a^{-,j_2}_{i_1+i_2,j_3;i_1,j_1}\eqno(3.8)$$
for $i_1,i_2\in\Bbb{Z}_2,\;j_1\in I_{i_1},\;j_2\in I_{i_2}$ and $j_3\in I_{i_1+i_2}$ by (2.32) and (3.3). Thus we have
$$ \vs^+_{i_1,j_1}\vs^-_{i_2,j_2}=\sum_{j_3\in I_{i_1+i_2}}((-1)^{i_2}a^{-,j_1}_{i_2,j_2;i_1+i_2,j_3}\vs^+_{i_1+i_2,j_3}+(-1)^{i_1}a^{+,j_2}_{i_1+i_2,j_3;i_1,j_1}\vs^-_{i_1+i_2,j_3}),\eqno(3.9)$$
$$\vs^-_{i_1,j_1}\vs^+_{i_2,j_2}=\sum_{j_3\in I_{i_1+i_2}}(a^{+,j_1}_{i_2,j_2;i_1+i_2,j_3}\vs^-_{i_1+i_2,j_3}+a^{-,j_2}_{i_1+i_2,j_3;i_1,j_1}\vs^+_{i_1+i_2,j_3})\eqno(3.10)$$
for $i_1,i_2\in\Bbb{Z}_2,\;j_1\in I_{i_1}$ and $j_2\in I_{i_2}$.

Let $i_p\in\Bbb{Z}_2,\;j_p\in I_{i_p}$ with $p=1,2,3.$ We have
\begin{eqnarray*}& &(\vs^+_{i_1,j_1}\vs^+_{i_2,j_2})\vs^-_{i_3,j_3}\\ &=&\sum_{j_4\in I_{i_1+i_2}}a^{+,j_4}_{i_1,j_1;i_2,j_2}\vs^+_{i_1+i_2,j_4}\vs^-_{i_3,j_3}\\ &=&\sum_{j_4\in I_{i_1+i_2},\;j_5\in I_{i_1+i_2+i_3}}a^{+,j_4}_{i_1,j_1;i_2,j_2}((-1)^{i_3}a^{-,j_4}_{i_3,j_3;i_1+i_2+i_3,j_5}\vs^+_{i_1+i_2+i_3,j_5}\\& &+(-1)^{i_1+i_2}a^{+,j_3}_{i_1+i_2+i_3,j_5;i_1+i_2,j_4}\vs^-_{i_1+i_2+i_3,j_5}),\hspace{7cm}(3.11)\end{eqnarray*}
\begin{eqnarray*}& &\vs^+_{i_1,j_1}(\vs^+_{i_2,j_2}\vs^-_{i_3,j_3})\\ &=&\sum_{j_4\in I_{i_2+i_3}}((-1)^{i_3}a^{-,j_2}_{i_3,j_3;i_2+i_3,j_4}\vs^+_{i_1,j_1}\vs^+_{i_2+i_3,j_4}+(-1)^{i_2}a^{+,j_3}_{i_2+i_3,j_4;i_2,j_2}\vs^+_{i_1,j_1}\vs^-_{i_2+i_3,j_4})\\ &=&\sum_{j_4\in I_{i_2+i_3},\;j_5\in I_{i_1+i_2+i_3}}[(-1)^{i_3}(a^{-,j_2}_{i_3,j_3;i_2+i_3,j_4}a^{+,j_5}_{i_1,j_1;i_2+i_3,j_4}\\ & &+a^{+,j_3}_{i_2+i_3,j_4;i_2,j_2}a^{-,j_1}_{i_2+i_3,j_4;i_1+i_2+i_3,j_5})\vs^+_{i_1+i_2+i_3,j_5}\\ & &+(-1)^{i_1+i_2}a^{+,j_3}_{i_2+i_3,j_4;i_2,j_2}
a^{+,j_4}_{i_1+i_2+i_3,j_5;i_1,j_1}\vs^-_{i_1+i_2+i_3,j_5}],\hspace{5.5cm}(3.12)\end{eqnarray*}
\begin{eqnarray*}& &(\vs^+_{i_1,j_1}\vs^-_{i_2,j_2})\vs^+_{i_3,j_3}\\ &=&\sum_{j_4\in I_{i_1+i_2}}((-1)^{i_2}a^{-,j_1}_{i_2,j_2;i_1+i_2,j_4}\vs^+_{i_1+i_2,j_4}+(-1)^{i_1}a^{+,j_2}_{i_1+i_2,j_4;i_1,j_1}\vs^-_{i_1+i_2,j_4})\vs^+_{i_3,j_3}\\&=&
\sum_{j_4\in I_{i_1+i_2},\;j_5\in I_{i_1+i_2+i_3}}[((-1)^{i_2}a^{+,j_5}_{i_1+i_2,j_4;i_3,j_3}a^{-,j_1}_{i_2,j_2;i_1+i_2,j_4}\\& &+(-1)^{i_1}a^{-,j_3}_{i_1+i_2+i_3,j_5;i_1+i_2,j_4}a^{+,j_2}_{i_1+i_2,j_4;i_1,j_1})\vs^+_{i_1+i_2+i_3,j_5}\\ && +(-1)^{i_1}a^{+,j_4}_{i_3,j_3;i_1+i_2+i_3,j_5}a^{+,j_2}_{i_1+i_2,j_4;i_1,j_1}\vs^-_{i_1+i_2+i_3,j_5}],\hspace{5.9cm}(3.13)\end{eqnarray*}
\begin{eqnarray*}& &\vs^+_{i_1,j_1}(\vs^-_{i_2,j_2}\vs^+_{i_3,j_3})\\ &=&\sum_{j_4\in I_{i_2+i_3}}\vs^+_{i_1,j_1}(a^{+,j_2}_{i_3,j_3;i_2+i_3,j_4}\vs^-_{i_2+i_3,j_4}+a^{-,j_3}_{i_2+i_3,j_4;i_2,j_2}\vs^+_{i_2+i_3,j_4})\\&=& \sum_{j_4\in I_{i_1+i_2},\;j_5\in I_{i_1+i_2+i_3}}[((-1)^{i_2+i_3}a^{-,j_1}_{i_2+i_3,j_4;i_1+i_2+i_3,j_5}a^{+,j_2}_{i_3,j_3;i_2+i_3,j_4}+a^{+,j_5}_{i_1,j_1;i_2+i_3,j_4}a^{-,j_3}_{i_2+i_3,j_4;i_2,j_2})\\ & &\times \vs^+_{i_1+i_2+i_3,j_5}+ (-1)^{i_1}a^{+,j_4}_{i_1+i_2+i_3,j_5;i_1,j_1}a^{+,j_2}_{i_3,j_3;i_2+i_3,j_4}\vs^-_{i_1+i_2+i_3,j_5}].\hspace{3.8cm}(3.14)\end{eqnarray*}
Thus we obtain
\begin{eqnarray*}& &\sum_{j_4\in I_{i_1+i_2}}a^{+,j_4}_{i_1,j_1;i_2,j_2}a^{-,j_4}_{i_3,j_3;i_1+i_2+i_3,j_5}\\ &=&
\sum_{j_4\in I_{i_2+i_3}}(a^{-,j_2}_{i_3,j_3;i_2+i_3,j_4}a^{+,j_5}_{i_1,j_1;i_2+i_3,j_4}+a^{+,j_3}_{i_2+i_3,j_4;i_2,j_2}a^{-,j_1}_{i_2+i_3,j_4;i_1+i_2+i_3,j_5})\hspace{2.8cm}(3.15)\end{eqnarray*}
by (3.11) and (3.12), and
\begin{eqnarray*}& &\sum_{j_4\in I_{i_1+i_2}}((-1)^{i_2}a^{+,j_5}_{i_1+i_2,j_4;i_3,j_3}a^{-,j_1}_{i_2,j_2;i_1+i_2,j_4}+(-1)^{i_1}a^{-,j_3}_{i_1+i_2+i_3,j_5;i_1+i_2,j_4}a^{+,j_2}_{i_1+i_2,j_4;i_1,j_1})\\ &=&\sum_{j_4\in I_{i_1+i_2}}((-1)^{i_2+i_3}a^{-,j_1}_{i_2+i_3,j_4;i_1+i_2+i_3,j_5}a^{+,j_2}_{i_3,j_3;i_2+i_3,j_4}+a^{+,j_5}_{i_1,j_1;i_2+i_3,j_4}a^{-,j_3}_{i_2+i_3,j_4;i_2,j_2})\hspace{1cm}(3.16)\end{eqnarray*}
by (3.13) and (3.14), for $i_p\in\Bbb{Z}_2,\;j_p\in I_{i_p}$ with $p=1,2,3$ and $j_5\in I_{i_1+i_2+i_3}$.
Conversely, we have the following conclusion.
\psp 

{\bf Proposition 3.1}. {\it Suppose that we have two} $\Bbb{Z}_2$-{\it graded associative algebra operations} $\circ_+$ {\it and} $\circ_-$ {\it on a} $\Bbb{Z}_2$-{\it graded vector space} ${\cal B}={\cal B}_0\oplus {\cal B}_1$ {\it such that under a basis} $\{\vt_{i,j}\mid i\in\Bbb{Z}_2,\;j\in I_i\}$ {\it of} ${\cal B}$, {\it the structure constants} 
$$\{a_{i_1,j_1;i_2,j_2}^{\pm,j_3}\mid i_1,i_2\in\Bbb{Z}_2,\;j_1\in I_{i_1},\;j_2\in I_{i_2},\;j_3\in I_{i_1+i_2}\}\eqno(3.17)$$
{\it satisfy (3.15) and (3.16), where}
$$\vt_{i_1,j_1}\circ_{\pm}\vt_{i_2,j_2}=\sum_{j_3\in I_{i_1+i_2}}a^{\pm,j_3}_{i_1,j_1;i_2,j_2}\vt_{i_1+i_2,j_3}\qquad\mbox{\it for}\;\;i_1,i_2\in\Bbb{Z}_2,\;j_1\in I_{i_1},\;j_2\in I_{i_2}.\eqno(3.18)$$
{\it Let} ${\cal A}^{\pm}_i$ {\it be the vector spaces with a basis} $\{\vs^{\pm}_{i,j}\mid j\in I_i\}$ {\it for} $i\in\Bbb{Z}_2$. {\it Set}
$${\cal A}={\cal A}_0\oplus {\cal A}_1={\cal A}^+\oplus {\cal A}^-\eqno(3.19)$$
{\it with}
$${\cal A}_i={\cal A}^+_i+{\cal A}^-_i,\;{\cal A}^{\pm}={\cal A}^{\pm}_0\oplus
{\cal A}^{\pm}_1.\eqno(3.20)$$
{\it We define the multiplication operation on} ${\cal A}$ {\it by (3.4), (3.9)and (3.10), and the bilinear form} $\la\cdot,\cdot\ra$ {\it by}
$$\la {\cal A}^{\pm},{\cal A}^{\pm}\ra=\{0\}\eqno(3.21)$$
{\it and}
$$\la \vs_{i_1,j_1}^-,\vs_{i_2,j_2}^+\ra=(-1)^{i_1}\la\vs_{i_2,j_2}^+, \vs_{i_1,j_1}^-\ra=
\dlt_{i_1,i_2}\dlt_{j_1,j_2}\eqno(3.22)$$
{\it for} $i_1,i_2\in\Bbb{Z}_2,\;j_1\in I_{i_1},\;j_2\in I_{i_2}$. {\it Then} ${\cal A}$ {\it forms a} $\Bbb{Z}_2$-{\it graded polarized associative algebra}.
\pse

{\it Proof}. By (3.5)-(3.8), (3.11)-(3.14) and the symmetry of (3.15) and (3.16) with respect to the signs ``$+$'' and ``$-$'', we only need to verify
$$(\vs^+_{i_1,j_1}\vs^-_{i_2,j_2})\vs^-_{i_3,j_3}=\vs^+_{i_1,j_1}(\vs^-_{i_2,j_2}\vs^-_{i_3,j_3}),\eqno(3.23)$$
$$\la \vs^+_{i_1,j_1}\vs^-_{i_2,j_2},\vs^+_{i_1+i_2,j_3}\ra=\la \vs^+_{i_1,j_1},\vs^-_{i_2,j_2}\vs^+_{i_1+i_2,j_3}\ra\eqno(3.24)$$
for $i_p\in\Bbb{Z}_2$ and $j_p\in I_p$. Note
\begin{eqnarray*}& &(\vs^+_{i_1,j_1}\vs^-_{i_2,j_2})\vs^-_{i_3,j_3}\\ &=&\sum_{j_4\in I_{i_1+i_2}}((-1)^{i_2}a^{-,j_1}_{i_2,j_2;i_1+i_2,j_4}\vs^+_{i_1+i_2,j_4}+(-1)^{i_1}a^{+,j_2}_{i_1+i_2,j_4;i_1,j_1}\vs^-_{i_1+i_2,j_4})\vs^-_{i_3,j_3}\\ &=&\sum_{j_4\in I_{i_1+i_2},\;j_5\in I_{i_1+i_2+i_3}}((-1)^{i_2+i_3}a^{-,j_1}_{i_2,j_2;i_1+i_2,j_4}a^{-,j_4}_{i_3,j_3;i_1+i_2+i_3,j_5}\vs^+_{i_1+i_2+i_3,j_5}\\& &+
(-1)^{i_1}(a^{-,j_1}_{i_2,j_2;i_1+i_2,j_4}a^{+,j_3}_{i_1+i_2+i_3,j_5;i_1+i_2,j_4}+a^{-,j_5}_{i_1+i_2,j_4;i_3,j_3}a^{+,j_2}_{i_1+i_2,j_4;i_1,j_1})\vs^-_{i_1+i_2+i_3,j_5},\hspace{0.8cm}(3.25)\end{eqnarray*}
\begin{eqnarray*}& &\vs^+_{i_1,j_1}(\vs^-_{i_2,j_2}\vs^-_{i_3,j_3})\\ &=&\sum_{j_4\in I_{i_2+i_3}}
a^{-,j_4}_{i_2,j_2;i_3,j_3}\vs^+_{i_1,j_1}\vs^-_{i_2+i_3,j_4}\\ &=& \sum_{j_4\in I_{i_2+i_3},\;j_5\in I_{i_1+i_2+i_3}}a^{-,j_4}_{i_2,j_2;i_3,j_3}((-1)^{i_2+i_3}a^{-,j_1}_{i_2+i_3,j_4;i_1+i_2+i_3,j_5}\vs^+_{i_1+i_2+i_3,j_5}\\& &+(-1)^{i_1}a^{+,j_4}_{i_1+i_2+i_3,j_5;i_1,j_1}\vs^-_{i_1+i_2+i_3,j_5}).\hspace{8cm}(3.26)\end{eqnarray*}
So (3.23) follows from (3.25), (3.26) and (3.15) with the change of indices:
$$ i_1\rta i_1+i_2+i_3\rta i_3\rta i_2\rta i_1,\;\;\;j_1\rta j_5\rta j_3\rta j_2\rta j_1.\eqno(3.27)$$
Moreover,
$$\la \vs^+_{i_1,j_1}\vs^-_{i_2,j_2},\vs^+_{i_1+i_2,j_3}\ra=(-1)^{i_1}a^{+,j_2}_{i_1+i_2,j_3;i_1,j_1}\eqno(3.28)$$
by (3.9) and (3.22), and 
$$\la \vs^+_{i_1,j_1},\vs^-_{i_2,j_2}\vs^+_{i_1+i_2,j_3}\ra=(-1)^{i_1}a^{+,j_2}_{i_1+i_2,j_3;i_1,j_1}
\eqno(3.29)$$
by (3.10) and (3.22).
Therefore, (3.24) holds.$\qquad\Box$
\psp

Let ${\cal A}$ be a polarized $\Bbb{Z}_2$-graded associative algebra with the bilinear form $\la\cdot,\cdot\ra_1$ and let ${\cal B}$ be a $\Bbb{Z}_2$-graded associative algebra with a nondegenerate supersymmetric associative bilinear form $\la\cdot,\cdot\ra_2$ (cf. (2.32)). Here ${\cal B}$ may not be polarized. Set
$$\td{\cal A}_0={\cal A}_0\otimes_{\Bbb{F}} {\cal B}_0+{\cal A}_1\otimes_{\Bbb{F}} {\cal B}_1,\;\;\td{\cal A}_1={\cal A}_0\otimes_{\Bbb{F}}{\cal B}_1+{\cal A}_1\otimes_{\Bbb{F}} {\cal B}_0,\eqno(3.30)$$
$$\td{\cal A}^+={\cal A}^+\otimes_{\Bbb{F}} {\cal B},\;\;\td{\cal A}^-={\cal A}^-\otimes_{\Bbb{F}} {\cal B}\eqno(3.31)$$
and
$$\td{\cal A}=\td{\cal A}_0\oplus \td{\cal A}_1=\td{\cal A}^+\oplus \td{\cal A}^-.\eqno(3.32)$$
Define the multiplication and bilinear form on $\td{\cal A}$ by
$$(a_1\otimes b_1)(a_2\otimes b_2)=a_1a_2\otimes b_1b_2,\;\;\la a_1\otimes b_1,a_2\otimes b_2\ra=(-1)^{i_1i_2}\la a_1,a_2\ra_1\la b_1,b_2\ra_2\eqno(3.33)$$
for $a_1\in{\cal A},\;a_2\in{\cal A}_{i_1}$ and $b_1\in{\cal B}_{i_2},\;b_2\in{\cal B}$. It is straightforward to verify the following proposition.
\psp

{\bf Proposition 3.2}. {\it The space} $\td{\cal A}$ {\it forms a polarized} $\Bbb{Z}_2$-{\it graded associative algebra}.
\psp

{\bf Example 3.1}. In the algebra $\mbb{F}[t,t^{-1}]$ of Laurent polynomials, we define the bilinear form
$$\la t^m,t^n\ra=\dlt_{m+n,-1}\qquad\for\;\;m,n\in\mbb{Z}.\eqno(3.34)$$
Set
$$(\mbb{F}[t,t^{-1}])^+=\mbb{F}[t],\qquad (\mbb{F}[t,t^{-1}])^-=\mbb{F}[t^{-1}]t^{-1}.\eqno(3.35)$$
Then $\mbb{F}[t,t^{-1}]=(\mbb{F}[t,t^{-1}])^+\oplus (\mbb{F}[t,t^{-1}])^-$ forms a polarized associative algebra (with $(\mbb{F}[t,t^{-1}])_0=\mbb{F}[t,t^{-1}]$ and $(\mbb{F}[t,t^{-1}])_1=\{0\}$).
\pse

{\bf Example 3.2}. Let $k$ be a positive integer and let $k_1\in\ol{0,k-1}$. Denote by $M_{k\times k}(\Bbb{F})$ the algebra of $k\times k$ matrices with entries in $\Bbb{F}$, and by $E_{j,l}$ the matrix with 1 as its $(j,l)$-entry and $0$ as the others. Define
$$M_{k\times k}(\Bbb{F})_0=\sum_{j,l\in\ol{1,k_1}}\Bbb{F}E_{j,l}+\sum_{p,q\in\ol{k_1+1,k}}\Bbb{F}E_{p,q},\eqno(3.36)$$
$$M_{k\times k}(\Bbb{F})_1=\sum_{j\in\ol{1,k_1},\;p\in\ol{k_1+1,k}}(\Bbb{F}E_{j,p}+\Bbb{F}E_{p,j})\eqno(3.37)$$
and
$$\tr\:A=\sum_{j\in\ol{1,k_1}}a_{j,j}-\sum_{p\in\ol{k_1+1,k}}a_{p,p}\qquad\for\;\;A=\sum_{j,l=1}^ka_{j,l}E_{j,l}\in M_{k\times k}(\Bbb{F}).\eqno(3.38)$$
Moreover, we define
$$\la A,B\ra=\tr\:AB\qquad\for\;\;A,B\in M_{k\times k}(\Bbb{F}).\eqno(3.39)$$
Then $M_{k\times k}(\Bbb{F})$ forms a $\Bbb{Z}_2$-graded associative algebra with the supersymmetric nondegenerate associative bilinear form $\la\cdot,\cdot\ra$.

According to Proposition 3.2, 
$$M_{k\times k}(\Bbb{F})\otimes_{\mbb{F}}\mbb{F}[t,t^{-1}]\cong M_{k\times k}(\Bbb{F}[t,t^{-1}])\eqno(3.40)$$
forms a polarized $\Bbb{Z}_2$-graded associative algebra, where $M_{k\times k}(\Bbb{F}[t,t^{-1}])$ denote the algebra of $k\times k$ matrices with entries in $\Bbb{F}[t,t^{-1}].$
\pse

{\bf Example 3.3}. Let $k>1$ be integer and $\Bbb{F}=\Bbb{C}$ the field of complex numbers. The {\it Hecke algebra} ${\cal H}_k$ is an associative algbera generated by $\{T_1,...,T_{k-1}\}$ with the following defining relations
$$T_iT_j=T_jT_i\qquad\mbox{whenever}\;\;|i-j|\geq 2,\eqno(3.41)$$
$$T_iT_{i+1}T_i=T_{i+1}T_iT_{i+1},\;\;T_i^2=(q-1)T_i+q\eqno(3.42)$$
for $i,j\in \ol{1,k-1}$, where $0\neq q\in \Bbb{C}$. Let $\zeta\in \Bbb{C}$ be a fixed constant.  According to Section 5 in [HKW], there exists a unique trace map ``$\tr$'' of ${\cal H}_k$ such that
$$\tr\:(1_{{\cal H}_k})=1,\;\;\tr\:(a T_nb)=\zeta \tr\:(ab)\qquad\mbox{for}\;\;a,b\in {\cal H}_n\eqno(3.43)$$
with $n\in \ol{1,k-2}$. This trace map is the key to define the well-known ``Jones polynomials'' of knots (e.g., cf. [HKW]). Furthermore, we define the bilinear form
$$\la u,v\ra=\tr\:uv\qquad\for\;\;u,v\in {\cal H}_k.\eqno(3.44)$$
Then $\la\cdot,\cdot\ra$ is a nondegenerate associative symmetric bilinear form  of ${\cal H}_k$ under a certain condition on $\zeta$. 

According to Proposition 3.2 and (3.40),
$$ {\cal H}_k\otimes_{\mbb{F}}M_{k\times k}(\Bbb{F}[t,t^{-1}])\eqno(3.45)$$
 forms a  $\Bbb{Z}_2$-graded polarized associative algebras under a certain condition of $\zeta$. Here we treat the odd part of ${\cal H}_k$ as zero.
\pse

{\bf Example 3.4}. Let $G$ be a group. Take a map $\ves: G\times G\rta \Bbb{F}^{\times}=\Bbb{F}\setminus\{0\}$ such that
$$\ves(g_1,g_2)\ves (g_1g_2,g_3)=\ves(g_1,g_2g_3)\ves(g_2,g_3)\qquad\for\;\;g_1,g_2,g_3\in G.\eqno(3.46)$$
Let $\Bbb{F}[G]_{\ves}$ be a vector space with a basis $\{u_g\mid g\in G\}$. Define the multiplication on $\Bbb{F}[G]_{\ves}$ by
$$u_{g_1}u_{g_2}=\ves(g_1,g_2)u_{g_1g_2}\qquad\for\;\;g_1,g_2\in G.\eqno(3.47)$$
Then $\Bbb{F}[G]_{\ves}$ forms an associative algebra, which is called a {\it twisted group algebra} of $G$. Moreover, we define $\tr: \Bbb{F}[G]_{\ves}\rta \Bbb{F}$ by
$$\tr\:(u_g)=\dlt_{g,0}\qquad\for\;\;g\in G,\eqno(3.48)$$
and
$$\la u,v\ra =\tr\:(uv)\qquad\for\;\;u,v\in \Bbb{F}[G]_{\ves}.\eqno(3.49)$$
It is straightforward to verify that $\la\cdot,\cdot\ra$ is a nondegenerate symmetric associative bilinear form of $\Bbb{F}[G]_{\ves}$.

According to Proposition 3.2 and (3.45), 
$$ \Bbb{F}[G]_{\ves}\otimes_{\mbb{F}}{\cal H}_k\otimes_{\mbb{F}}M_{k\times k}(\Bbb{F}[t,t^{-1}])\eqno(3.50)$$
 forms a  $\Bbb{Z}_2$-graded polarized associative algebra. Here we treat the odd part of $\Bbb{F}[G]_{\ves}$ as zero.
\psp

The following work is necessary for defining  Poisson superpairs over an infinite-dimensional $\Bbb{Z}_2$-graded polarized associative algebra.
 
Let ${\cal A}$ be a $\Bbb{Z}_2$-graded polarized associative algebra with the bilinear form $\la\cdot,\cdot\ra$. Define
$${\cal A}^{(0)}=\{u\in{\cal A}^+\mid {\cal A}^-u\subset {\cal A}^-\}\eqno(3.51)$$
and
$${\cal A}^{(m+1)}=\{u\in{\cal A}^+\mid {\cal A}^-u\subset {\cal A}^{(m)}+{\cal A}^-\}\eqno(3.52)$$
for $m\in\Bbb{N}$ by induction. Then we have
$${\cal A}^{(n)}\subset {\cal A}^{(n+1)}\qquad\for\;\;n\in\Bbb{N}.\eqno(3.53)$$
Set
$${\cal A}^{(-1)}={\cal A}^-,\;\;{\cal A}^{(-2-m)}=\{u\in {\cal A}^-\mid \la u,{\cal A}^{(m)}\ra=\{0\}\}\qquad\for\;\;m\in\Bbb{N}.\eqno(3.54)$$
Then (3.53) also holds for a negative intger $n<-1$.
Since ${\cal A}^+$ and ${\cal A}^-$ are $\Bbb{Z}_2$-graded subalgebras of ${\cal A}$, all ${\cal A}^{(n)}$ with $n\in\Bbb{Z}$ are $\Bbb{Z}_2$-graded subspaces of ${\cal A}$ by (2.32).
\psp

{\bf Proposition 3.3}. {\it For} $m,n\in\Bbb{N}$, 
$${\cal A}^{(m)}{\cal A}^{(n)}\subset {\cal A}^{(m+n)},\eqno(3.55)$$
$${\cal A}^{(m+n)}{\cal A}^{(-n-1)},\:{\cal A}^{(-n-1)}{\cal
A}^{(m+n)}\subset {\cal A}^{(m-1)}+{\cal A}^-.\eqno(3.56)$$
$${\cal A}^{(m)}{\cal A}^{(-m-n-2)},\:{\cal A}^{(-m-n-2)}{\cal
A}^{(m)}\subset {\cal A}^{(-n-2)},\eqno(3.57)$$
$${\cal A}^{(-m-1)}{\cal A}^{(-n-1)}\subset {\cal A}^{(-m-n-2)}.\eqno(3.58)$$

{\it Proof}. We prove (3.55) by induction on $m$. Note
$${\cal A}^-({\cal A}^{(0)}{\cal A}^{(n)})=({\cal A}^-{\cal A}^{(0)}){\cal A}^{(n)}\subset {\cal A}^-{\cal A}^{(n)}\subset {\cal A}^{(n-1)}+{\cal A}^-\eqno(3.59)$$
by (3.51) and (3.52). Moreover, (3.52) and (3.59) imply (3.55) with $m=0$. Suppose that (3.55) holds for $m=k$ with $k\in\Bbb{N}$. We have
$${\cal A}^-({\cal A}^{(k+1)}{\cal A}^{(n)})=({\cal A}^-{\cal A}^{(k+1)}){\cal A}^{(n)}\subset ({\cal A}^{(k)}+{\cal A}^-){\cal A}^{(n)}\subset {\cal A}^{(k+n)}+{\cal A}^-\eqno(3.60)$$
by (3.52). Again (3.52) and (3.60) imply (3.55) with $m=k+1$. So (3.55) holds.

We define
$${\cal A}^{[0]}=\{u\in {\cal A}^+\mid u{\cal A}^- \subset {\cal A}^-\}\eqno(3.61)$$
and
$${\cal A}^{[m+1]}=\{u\in{\cal A}^+\mid u{\cal A}^-\subset {\cal A}^{[m]}+{\cal A}^-\}\eqno(3.62)$$
for $m\in\Bbb{N}$ by induction. Note
$$\la {\cal A}^-, {\cal A}^{(0)}{\cal A}^-\ra= \la {\cal A}^-{\cal A}^{(0)},{\cal A}^-\ra\subset\la {\cal A}^-,{\cal A}^-\ra=\{0\}\eqno(3.63)$$
by (2.32) and (2.33). By the nondegeneracy of $\la\cdot,\cdot\ra$ and (2.34), we have
$${\cal A}^{(0)}{\cal A}^-\subset {\cal A}^-.\eqno(3.64)$$
Thus 
$${\cal A}^{(0)}\subset {\cal A}^{[0]}.\eqno(3.65)$$
By the symmetric proof as that in the above, we also have ${\cal A}^{[0]}\subset {\cal A}^{(0)}$. Hence
$${\cal A}^{(0)}={\cal A}^{[0]}.\eqno(3.66)$$

Suppose that
$${\cal A}^{(m)}={\cal A}^{[m]}\eqno(3.67)$$
for  $m \leq k$ with $k\in\Bbb{N}$. We have
\begin{eqnarray*} & &{\cal A}^-({\cal A}^{(k+1)}{\cal A}^-)=({\cal A}^-{\cal A}^{(k+1)}){\cal A}^-\subset({\cal A}^{(k)}+{\cal A}^-){\cal A}^-\\& &
\subset {\cal A}^{[k]}{\cal A}^-+{\cal A}^-\subset {\cal A}^{[k-1]}+{\cal A}^-={\cal A}^{(k-1)}+{\cal A}^-\hspace{5.9cm}(3.68)\end{eqnarray*}
by (3.52), (3.62) and (3.67). Definition (3.52) implies
$${\cal A}^{(k+1)}{\cal A}^-\subset {\cal A}^{(k)}={\cal A}^{[k]}.\eqno(3.69)$$
By (3.62), we have
$${\cal A}^{(k+1)}\subset {\cal A}^{[k+1]}.\eqno(3.70)$$
A symmetric argument shows ${\cal A}^{[k+1]}\subset {\cal A}^{(k+1)}$. Hence
$${\cal A}^{(k+1)}={\cal A}^{[k+1]}.\eqno(3.71)$$
Therefore, (3.67) holds for any $m\in\Bbb{N}$ by induction. 

Expressions (3.52), (3.62) and (3.67) imply that (3.56) holds for $n=0$. Assume $n>0$. Then we have
$$\la {\cal A}^-, {\cal A}^{(n)}{\cal A}^{(-n-1)}\ra\subset \la {\cal A}^-{\cal A}^{(n)}, {\cal A}^{(-n-1)}\ra\subset \la {\cal A}^{(n-1)}+{\cal A}^-, {\cal A}^{(-n-1)}\ra=\{0\}\eqno(3.72)$$
by (2.32), (2.33) and (3.54). The nondegeneracy of $\la\cdot,\cdot\ra$ and (2.34) imply
$${\cal A}^{(n)}{\cal A}^{(-n-1)}\subset{\cal A}^-= {\cal A}^{(-1)}.\eqno(3.73)$$
Suppose that
$${\cal A}^{(m+n)}{\cal A}^{(-n-1)}\subset {\cal A}^{(m-1)}+{\cal A}^-\eqno(3.74)$$
holds for some $m\in\Bbb{N}$. We have
\begin{eqnarray*}& &{\cal A}^-({\cal A}^{(m+1+n)}{\cal A}^{(-n-1)})=({\cal A}^-{\cal A}^{(m+1+n)}){\cal A}^{(-n-1)}\subset({\cal A}^{(m+n)}+{\cal A}^-){\cal A}^{(-n-1)}\\& & \subset{\cal A}^{(m+n)}{\cal A}^{(-n-1)}+{\cal A}^-
 \subset{\cal A}^{(m-1)}+{\cal A}^-\hspace{7.1cm}(3.75)\end{eqnarray*}
by (3.52), (3.74) and the fact ${\cal A}^-$ is a subalgebra of ${\cal A}$ (note ${\cal A}^{(-n-1)}\subset {\cal A}^-$ by (3.54)). Again the definition (3.52)  imply
$${\cal A}^{(m+1+n)}{\cal A}^{(-n-1)}\subset  {\cal A}^{(m)}+{\cal A}^-.\eqno(3.76)$$
By induction, (3.74) holds for any $m\in\Bbb{N}$. By means of (3.67), we can symmetrically prove
$${\cal A}^{(-n-1)}{\cal A}^{(m+n)}\subset {\cal A}^{(m-1)}+{\cal A}^-\eqno(3.77)$$
This proves (3.56).

For $m,n\in\Bbb{N}$, we have 
$$\la {\cal A}^{(n)},{\cal A}^{(m)}{\cal A}^{(-m-n-2)}\ra=\la {\cal A}^{(n)}{\cal A}^{(m)},{\cal A}^{(-m-n-2)}\ra\subset \la {\cal A}^{(m+n)},{\cal A}^{(-m-n-2)}\ra=\{0\}\eqno(3.78)$$
by (2.32), (3.54) and (3.55). Moreover, (3.54) and (3.78) imply
$${\cal A}^{(m)}{\cal A}^{(-m-n-2)}\subset {\cal A}^{(-n-2)}.\eqno(3.79)$$
Symmetrically, we can prove
$${\cal A}^{(-m-n-2)}{\cal A}^{(m)}\subset {\cal A}^{(-n-2)}\eqno(3.80)$$
by (3.67). So (3.57) holds.

Observe that
\begin{eqnarray*}\hspace{2cm} & &\la{\cal A}^{(m+n)}, {\cal A}^{(-m-1)}{\cal A}^{(-n-1)} \ra=\la{\cal A}^{(m+n)} {\cal A}^{(-m-1)},{\cal A}^{(-n-1)}\\& &\subset \la{\cal A}^{(n-1)} +{\cal A}^-,{\cal A}^{(-n-1)}\ra=\{0\}\hspace{6.4cm}(3.81)\end{eqnarray*}
by (2.32), (2.33), (3.54) and (3.56). Moreover, (3.54) and (3.81) imply
$${\cal A}^{(-m-1)}{\cal A}^{(-n-1)} \subset {\cal A}^{(-m-n-2)},\eqno(3.82)$$
that is, (3.58) holds.$\qquad\Box$

\section{Poisson Superpairs}

In this section, we shall construct Poisson superpairs over a $\Bbb{Z}_2$-graded polarized associative algebra.

Let $\Lmd$ be a vector space, not necessarily finite-dimensional. Let $F(\Lmd)$ be the free associative algebra generated by $\Lmd$. Then the exterior algebra ${\cal E}$ generated by $\Lmd$ is isomorphic to 
$$ {\cal E}=F(\Lmd)/(\{uv+vu\mid u,v\in\Lmd\}).\eqno(4.1)$$
We can identify $\Lmd$ with its image in ${\cal E}$. Note that 
$${\cal E}={\cal E}_0\oplus  {\cal E}_1,\qquad\mbox{where}\;\;  {\cal E}_0=\sum_{n=0}^{\infty}\Lmd^{2n},\;\;{\cal E}_1=\sum_{n=0}^{\infty}\Lmd^{2n+1}.\eqno(4.2)$$
With respect to the above grading, ${\cal E}$ becomes a super commutative associative algebra, that is,
$$uv=(-1)^{ij} vu\qquad\for\;\;u\in{\cal E}_i,\;v\in{\cal E}_j.\eqno(4.3)$$

Let 
$${\cal A}={\cal A}^+\oplus {\cal A}^-={\cal A}_0\oplus{\cal A}_1\eqno(4.4)$$
be infinite-dimensional $\Bbb{Z}_2$-graded polarized associative algebra with the bilinear form $\la\cdot,\cdot\ra$. Recall the notations of ${\cal A}^{(n)}$ with $n\in\Bbb{Z}$ defined in (3.51), (3.52) and (3.54). Their properties, which are important to the following construction, have been presented in Proposition 3.3 (cf. (3.55)-(3.58)). Assume that the algebra ${\cal A}$ satisfies the following condition:
$${\cal A}^+=\bigcup_{m=0}^{\infty}{\cal A}^{(m)},\;\;\dim {\cal A}^{(m)}<\infty.\eqno(4.5)$$
We call such an algebra ${\cal A}$ a {\it semi-finitely-filtered polarized} $\Bbb{Z}_2$-{\it graded associative algebra}. Examples of this type of algebras have been given in Examples 3.1-2.4.

Set
$${\cal A}^{(m)}_i={\cal A}^{(m)}\bigcap {\cal A}_i,\;\;k_{i,m}=\dim {\cal A}^{(m)}_i\eqno(4.6)$$
for $i\in\Bbb{Z}_2$ and $m\in\Bbb{N}$. 
Take a basis $\{\vs_{i,j}\mid j\in J_i\}$ of ${\cal A}^+_i$ for $i\in\Bbb{Z}_2$ with  
$$J_i=\Bbb{N}+1\;\;\mbox{or}\;\;\ol{1,n}\;\;\mbox{for some}\;\;n\in\Bbb{N}+1\eqno(4.7)$$
(which is guaranteed by (4.5)) such that
$$\{\vs_{i,j}\mid j\in\ol{1,k_{i,m}}\}\;\;\mbox{is a basis of}\;\;{\cal A}^{(m)}_i\eqno(4.8)$$
for $i\in\Bbb{Z}_2$ and $m\in\Bbb{N}$.

Take a subset $\{\vs_{i,-j}\mid j\in J_i\}$ of ${\cal A}^-_i$ for $i\in\Bbb{Z}_2$ such that
$$\la \vs_{i_1,-j_1},\vs_{i_2,j_2}\ra=\dlt_{i_1,i_2}\dlt_{j_1,j_2}\qquad\for\;\;i_1,i_2\in\Bbb{Z}_2,\;j_1\in J_{i_1},\;j_2\in J_{i_2}\eqno(4.9)$$
and for any $u\in {\cal A}^-_i$, 
$$u=\sum_{i\in{Z}_2,\;j\in J_i}\lmd_{i,j}\vs_{i,-j}\qquad\mbox{with}\;\;\lmd_{i,j}\in\Bbb{F}.\eqno(4.10)$$
So $\{\vs_{i,-j}\mid j\in J_i\}$  is a {\it dual basis}  in ${\cal A}^-_i$ of $\{\vs_{i,j}\mid j\in J_i\}$, in the possible sense of toplogical completion. 
By (3.54)
$$\vs_{i,-j}\in {\cal A}^{(-m-2)}\qquad\for\;\;i\in\Bbb{Z}_2,\;k_{i,m}<j\in J_i.\eqno(4.11)$$

Pick a positive integer $\iota$. Set
$$I_i=(-J_i)\bigcap \ol{1,k_{i,\iota}}\qquad\for\;\;i\in\Bbb{Z}_2.\eqno(4.12)$$
For $i\in\Bbb{Z}_2$, let $\{x_{i,j}\mid j\in I_i\}$
be the variables taking values in ${\cal E}_i$ (cf. (4.2)). Denote by ${\cal P}$ the algebra of supersymmetric polynomials in $\{x_{i,j}\mid i\in\Bbb{Z}_2,\; j\in I_i\}$. Then ${\cal P}$ has the $\Bbb{Z}_2$-grading
$${\cal P}_i=\mbox{Span}\:\{x_{i_1,j_1}\cdots x_{i_\ell,j_\ell}\mid \ell\in\Bbb{N},\;i_p\in\Bbb{Z}_2,\;j_p\in I_{i_p},\;\sum_{r=1}^\ell i_r\equiv i\}\eqno(4.13)$$
for $i\in\Bbb{Z}_2$, and
$$fg=(-1)^{i_1i_2}fg\qquad\for\;\;f\in{\cal P}_{i_1},\;g\in{\cal P}_{i_2}.\eqno(4.14)$$

For $i_1\in\Bbb{Z}_2$ and $j_1\in I_{i_1}$, we define a linear transformation $\ptl_{i_1,j_1}$ on ${\cal P}$ by 
$$\ptl_{i_1,j_1}(fg)=\ptl_{i_1,j_1}(f)g+(-1)^{i_1i_2}f\ptl_{i_1,j_1}(g)\qquad\for\;\;i_2\in\Bbb{Z}_2,\;f\in{\cal P}_{i_2},\;g\in{\cal P}\eqno(4.15)$$ 
and
$$\ptl_{i_1,j_1}(x_{i_2,j_2})=\dlt_{i_1,i_2}\dlt_{j_1,j_2}\qquad\;\;\for\;\;i_2\in\Bbb{Z}_2,\;j_{i_2}\in I_{i_2}.\eqno(4.16)$$
Then $\ptl_{i_1,j_1}$ is a supersymmetric derivation of ${\cal P}$ with parity
$i_1$. Set
$${\cal W}_0=\{\sum_{i\in{Z}_2,\;j\in I_i}f_{i,j}\ptl_{i,j}\mid f_{i,j}\in {\cal P}_i\},\;\;{\cal W}_1=\{\sum_{i\in{Z}_2,\;j\in I_i}f_{i,j}\ptl_{i,j}\mid f_{i,j}\in {\cal P}_{i+1}\}\eqno(4.17)$$
and
$${\cal W}={\cal W}_0+{\cal W}_1\eqno(4.18)$$
as a subspace of supersymmetric derivations of ${\cal P}$. The Lie superbracket on ${\cal W}$ is defined by
$$[d_1,d_2]=d_1d_2-(-1)^{i_1i_2}d_2d_1\qquad \for\;\;d_1\in {\cal W}_{i_1},\;d_2\in{\cal W}_{i_2}.\eqno(4.19)$$
In fact, if $d_1=\sum_{i\in{Z}_2,\;j\in I_i}f_{i,j}\ptl_{i,j}\in{\cal W}_{i_1}$ and $d_2=\sum_{i\in{Z}_2,\;j\in I_i}g_{i,j}\ptl_{i,j}\in{\cal W}_{i_2}$, then we have
$$[d_1,d_2]=\sum_{i\in{Z}_2,\;j\in I_i}\sum_{i_1\in{Z}_2,\;j_1\in I_{i_1}}[f_{i_1,j_1}\ptl_{i_1,j_1}(g_{i,j})-(-1)^{i_1i_2}g_{i_1,j_1}\ptl_{i_1,j_1}(f_{i,j})]\ptl_{i,j}.\eqno(4.20)$$

Define the vector space
$$\bar{\cal G}=\{\sum_{i\in{Z}_2,\;j\in J_i}\xi_{i,-j}\vs_{i,-j}\mid\xi_{i,-j}\in {\cal P}\}+\sum_{i\in{Z}_2,\;j\in J_i}{\cal P}\vs_{i,j}.\eqno(4.21)$$
Equip $\bar{\cal G}$ with the $\Bbb{Z}_2$-grading:
$$\bar{\cal G}_0=\{\sum_{i\in{Z}_2,\;j\in J_i}\xi_{i,-j}\vs_{i,-j}\mid\xi_{i,-j}\in {\cal P}_i\}+\sum_{i\in{Z}_2,\;j\in J_i}{\cal P}_i\vs_{i,j},\eqno(4.22)$$
$$\bar{\cal G}_1=\{\sum_{i\in{Z}_2,\;j\in J_i}\xi_{i,-j}\vs_{i,-j}\mid\xi_{i,-j}\in {\cal P}_{i+1}\}+\sum_{i\in{Z}_2,\;j\in J_i}{\cal P}_{i+1}\vs_{i,j}.\eqno(4.23)$$
Set
$$\bar{J}_i=J_i\bigcup(-J_i)\qquad\for\;\;i\in\Bbb{Z}_2.\eqno(4.24)$$
For covenience, we use the notation
$$\vs_{_{\xi}}=\sum_{i\in{Z}_2,\;j\in \bar{J}_i}\xi_{i,j}\vs_{i,j}\in\bar{\cal G}.\eqno(4.25)$$
Moreover, we define  the multiplication on $\bar{\cal G}$ by
$$\vs_{_{\xi}}\vs_{\eta}=\sum_{i_1,i_2\in{Z}_2,\;j_1\in \bar{J}_{i_1},\;j_2\in \bar{J}_{i_2}}(-1)^{i_1(i_2+p)}\xi_{i_1,j_1}\eta_{i_2,j_2}\vs_{i_1,j_1}\vs_{i_2,j_2}\eqno(4.26)$$
for $\vs_{_{\xi}}\in \bar{\cal G}$ and $\vs_{\eta}\in\bar{\cal G}_p$. The above expression is well defined because of Proposition 3.3. It can be verified that the space $\bar{\cal G}$ forms a $\Bbb{Z}_2$-graded associative algebra with respect to the multiplication in (4.26). Furthermore, we define a bilinear form $\bar{\cal G}$ by
$$\la\vs_{_{\xi}},\vs_{\eta}\ra=\sum_{i\in{Z}_2,\;j\in J_i}(-1)^{i(i+p)}(\xi_{i,-j}\eta_{i,j}+(-1)^i\xi_{i,j}\eta_{i,-j})\eqno(4.27)$$
for $\vs_{_{\xi}}\in \bar{\cal G}$ and $\vs_{\eta}\in\bar{\cal G}_p$, where the sum is finite by (4.21). 
It is straightforward to verify that the above bilinear form is a $\Bbb{Z}_2$-graded supersymmetric associative  bilinear form of $\bar{\cal G}$. In fact, one can view $\bar{\cal G}$ as an extension algebra of ${\cal A}$ with extended $\Bbb{Z}_2$-graded supersymmetric associative  bilinear form $\la\cdot,\cdot\ra$.

Define 
$${\cal G}=\{\sum_{i\in{Z}_2,\;j\in I_i}\xi_{i,j}\vs_{i,j}\mid\xi_{i,j}\in {\cal P}\}\eqno(4.28)$$
(cf. (4.12)). Then ${\cal G}$ forms a $\Bbb{Z}_2$-graded subspace of $\bar{\cal G}$, that is,
$${\cal G}={\cal G}_0\oplus{\cal G}_1,\qquad {\cal G}_i={\cal G}\bigcap\bar{\cal G}_i.\eqno(4.29)$$
In general,  ${\cal G}$ does not forms a subalgebra of $\bar{\cal G}$. We shall use the convention that
$$ \vs_{_{\xi}}\in {\cal G}\;\;\mbox{implies}\;\;\xi_{i,j}=0\;\;\for\;\;j>k_{i,\iota}.\eqno(4.30)$$
Thus 
$$\vs_{_{\xi}}=\sum_{i\in{Z}_2,\;j\in I_i}\xi_{i,j}\vs_{i,j}\qquad\mbox{if}\;\;\vs_{_{\xi}}\in{\cal G}\eqno(4.31)$$
(cf. (4.12)).

We define
$$\ptl_{\vs_{_{\xi}}}=\sum_{i\in{Z}_2,\;j\in I_i}\xi_{i,j}\ptl_{i,j}\qquad\for\;\;\vs_{_{\xi}}\in{\cal G}.\eqno(4.32)$$
The map
$$u\rta \ptl_u\eqno(4.33)$$
gives a $\Bbb{Z}_2$-graded linear isomorphism between ${\cal G}$ and ${\cal W}$. Moreover, we define the action of ${\cal W}$ on ${\cal G}$ by
$$d(\vs_{_{\xi}})=\sum_{i\in{Z}_2,\;j\in I_i}d(\xi_{i,j})\vs_{i,j}\qquad\for\;\;d\in{\cal W},\;\;\vs_{_{\xi}}\in{\cal G}.\eqno(4.34)$$

We can use (4.33) as an identification of ${\cal G}$ with ${\cal W}$ (cf. (4.18) and (4.20)). We denote by $[\cdot,\cdot]_0$ the corresponding Lie superbracket of ${\cal G}$. By 
(4.19), (4.20) and (4.33), we have
$$[u,v]_0=\ptl_u(v)-(-1)^{i_1i_2}\ptl_v(u)\qquad\for\;\;u\in{\cal G}_{i_1},\;v\in{\cal G}_{i_2}.\eqno(4.35)$$
The algebra ${\cal P}$ becomes a ${\cal G}$-module with the action
$$u(f)=\ptl_u(f)\qquad \for\;\; u\in{\cal G},\;f\in{\cal P}.\eqno(4.36)$$

Set
$$\Omega=\sum_{i\in{Z}_2,\;j\in I_i}{\cal P}\vs_{i,-j}.\eqno(4.37)$$
Then $\Omega$ forms a $\Bbb{Z}_2$-graded subspace of $\bar{\cal G}$, that is,
$$\Omega=\Omega_0\oplus\Omega_1,\qquad \Omega_i={\cal G}\bigcap{\cal G}_i.\eqno(4.38)$$
We shall also use the convention that
$$ \vs_{_{\xi}}\in \Omega\;\;\mbox{implies}\;\;\xi_{i,j}=0\;\;\for\;\;j<-k_{i,\iota}.\eqno(4.39)$$
Thus 
$$\vs_{_{\xi}}=\sum_{i\in{Z}_2,\;j\in I_i}\xi_{i,-j}\vs_{i,-j}\qquad\mbox{if}\;\;\vs_{_{\xi}}\in\Omega,\eqno(4.40)$$
where the sum is finite by (4.37). We identify $\Omega$ with a subspace of one-forms by
$$w(u)=\la u,w\ra \qquad\for\;\;w\in\Omega,\;u\in {\cal G}\eqno(4.41)$$
(cf. (4.27)). For $f=f_0+f_1$ with $f_0\in{\cal P}_0$ and $f_1\in{\cal P}_1$, we define
\begin{eqnarray*}\vs_{_{(f)}}&=&\sum_{i\in{Z}_2,\;j\in J_i}((-1)^i\ptl_{i,-j}(f_0)+\ptl_{i,-j}(f_1))\vs_{i,j}\\& &+\sum_{i\in{Z}_2,\;j\in\ol{1, k_{i,\iota}}}(\ptl_{i,j}(f_0)+(-1)^i\ptl_{i,j}(f_1))\vs_{i,-j}\in\Omega.\hspace{5.1cm}(4.42)\end{eqnarray*}
Then
$$df(u)=\ptl_u(f)=\la u,\vs_{_{(f)}}\ra\qquad\for\;\;u\in{\cal G},\;f\in{\cal P}.\eqno(4.43)$$	
Hence
$$df=\vs_{_{(f)}}\qquad \for\;\;f\in{\cal P}.\eqno(4.44)$$

Let $H:\Omega\rta {\cal G}$ be a map of the form
$$H(\vs_{\eta})=\sum_{i\in\Bbb{Z}_2,\;j\in I_i}\eta_{i,j}a_{i,j}\eqno(4.45)$$
with
$$a_{i,j}\in {\cal G}_i\eqno(4.46)$$
for $\vs_{\eta}\in\Omega$. Given $u\in{\cal G}_q$, we define the  map $\ptl_u(H):\Omega\rta {\cal G}$ by
$$\ptl_u(H)(\vs_{\eta})=\sum_{i_1,i_2\in\Bbb{Z}_2,\;j_1\in I_{i_1},\;j_2\in I_{i_2}}(-1)^{(i_2+p)q}\eta_{i_2,j_2}\ptl_u(a_{i_2,j_2})\eqno(4.47)$$
for $\vs_{\eta}\in\Omega_p$. Note that for $\vs_{\eta}\in\Omega_p$ and $u\in{\cal G}_q$, we have
$$\ptl_u(H(\vs_{\eta}))=\ptl_u(H)(\vs_{\eta})+H(\ptl_u(\vs_{\eta}))\eqno(4.48)$$
by (4.32), (4.34) and (4.47).
Define
$$\{f,g\}_H=\la H(\vs_{_{(f)}}),\vs_{_{(g)}}\ra\qquad\for\;\;f,g\in{\cal P}.\eqno(4.49)$$
By (3.59)-(3.71) in [X3] (also see Section 4 in [X1]) and (4.48), we have the following Lemma.
\psp

{\bf Lemma 4.1}. {\it The map} $\{\cdot,\cdot\}_H$ {\it forms a Poisson superbracket if} $H$ {\it satisfies}
$$\la H(v),u\ra=-(-1)^{i_1i_2}\la H(u),v\ra\eqno(4.50)$$
{\it and}
$$\la\ptl_{H(u)}(H)(v),w\ra+(-1)^{i_1(i_2+i_3)}\la\ptl_{H(v)}(H)(w),u\ra+(-1)^{(i_1+i_2)i_3}\la\ptl_{H(w)}(H)(u),v\ra=0\eqno(4.51)$$
{\it for} $u\in\Omega_{i_1},\;v\in\Omega_{i_2}$ {\it and} $w\in\Omega_{i_3}$. 
\psp

Denote
$$(\vs_{_{\xi}})_{\pm}=\sum_{i\in{Z}_2,\;j\in J_i}\xi_{i,\pm j}\vs_{i,\pm j}\qquad\for\;\;\vs_{_{\xi}}\in\bar{\cal G}.\eqno(4.52)$$
Set
$$\bar{\cal G}^{\pm}=\{u_{\pm}\mid u\in\bar{\cal G}\}.\eqno(4.53)$$
Then $\bar{\cal G}^+$ and $\bar{\cal G}^-$ form $\Bbb{Z}_2$-graded associative subalgebras of $\bar{\cal G}$. In fact,
$$\bar{\cal G}=\bar{\cal G}^+\oplus\bar{\cal G}^-\;\;\mbox{is a polarization of}\;\;\bar{\cal G}\eqno(4.54)$$
 with respect to the multiplication in (4.26) and the bilinear form in (4.27). For any $v\in\bar{\cal G}$, we write
$$v=v_++v_-\qquad\mbox{with}\;\;v_{\pm}\in\bar{\cal G}^{\pm}.\eqno(4.55)$$

In order to prove our main theorem in this section, we need the following lemma.
\psp

{\bf Lemma 4.2}. {\it For} $u\in\bar{\cal G}_{i_1},\;v\in\bar{\cal G}_{i_2}$ {\it and} $w\in\bar{\cal G}_{i_3}$, {\it we have}:
\begin{eqnarray*}\la u,vw\ra&=&\la u,v_+w_-\ra+(-1)^{i_1(i_2+i_3)}\la v,w_+u_-\ra+(-1)^{(i_1+i_2)i_3}\la w,u_+v_-\ra\\&=&\la u,v_-w^+\ra+(-1)^{i_1(i_2+i_3)}\la v,w_-u_+\ra+(-1)^{(i_1+i_2)i_3}\la w,u_-v_+\ra.\hspace{2.1cm}(4.56)\end{eqnarray*}

{\it Proof}. For an expression $h(u,v,w)$, we denote
$$h(u,v,w)+c.p.=h(u,v,w)+(-1)^{i_1(i_2+i_3)}h(v,w,u)+(-1)^{i_3(i_1+i_2)}h(w,u,v).\eqno(4.57)$$
Moreover, 
$$\la uv,w\ra=(-1)^{i_3(i_1+i_2)}\la w,uv\ra=(-1)^{i_3(i_1+i_2)}\la wu,v\ra=(-1)^{i_1(i_2+i_3)}\la v,wu\ra\eqno(4.58)$$
by the supersymmetry and associativity of $\la\cdot,\cdot\ra$ on $\bar{\cal G}$.

Note by (2.33), (2.34), (4.58), and the supersymmetry and associativity of $\la\cdot,\cdot\ra$ on $\bar{\cal G}$, 
\begin{eqnarray*}\la u,v_+w_-\ra+c.p.&=& \la u_++u_-,v_+w_-\ra+c.p.\\ &=&
\la u_+,v_+w_-\ra+\la u_-,v_+w_-\ra+c.p.\\&=& \la u_+v_+,w_-\ra+(-1)^{i_3(i_1+i_2)}\la w_-u_-,v_+\ra+c.p.\\&=& \la u_+v_+,w\ra+(-1)^{i_3(i_1+i_2)}\la w_-u_-,v\ra+c.p.\\ &=& \la u_+v_+,w\ra+\la u_-v_-,w\ra+c.p.\hspace{5.8cm}(4.59)\end{eqnarray*}
and
\begin{eqnarray*}\la u,v_-w_+\ra+c.p.&=& \la u_++u_-,v_-w_+\ra+c.p.\\ &=&\la u_+,v_-w_+\ra+\la u_-,v_-w_+\ra+c.p.\\ &=&
(-1)^{i_1(i_2+i_3)}\la v_-,w_+u_+\ra+\la u_-v_-,w_+\ra+c.p.\\ &=&
(-1)^{i_1(i_2+i_3)}\la v,w_+u_+\ra+\la u_-v_-,w\ra+c.p.\\ &=& \la u_+v_+,w\ra+\la u_-v_-,w\ra+c.p.\hspace{5.8cm}(4.60)\end{eqnarray*}
Thus 
\begin{eqnarray*}& &\la u,v_+w_-\ra+(-1)^{i_1(i_2+i_3)}\la v,w_+u_-\ra+(-1)^{(i_1+i_2)i_3}\la w,u_+v_-\ra\\&=&\la u,v_-w_+\ra+(-1)^{i_1(i_2+i_3)}\la v,w_-u_+\ra+(-1)^{(i_1+i_2)i_3}\la w,u_-v_+\ra\\&= &{1\over 3}[\la u,v_+w_-\ra+\la u,v_-w_+\ra+\la u_+v_+,w\ra+\la u_-v_-,w\ra+c.p.]\\&=& {1\over 3}[\la u,v_+w_-\ra+\la u,v_-w_+\ra+\la u,v_+w_+\ra+\la u,v_-w_-\ra+c.p.]\\&=& {1\over 3}[\la u,vw\ra
+c.p.]\\ &=& \la u,vw\ra\hspace{12.6cm}(4.61)\end{eqnarray*}
by the supersymmetry and associativity of $\la\cdot,\cdot\ra$ on $\bar{\cal G}.\qquad\Box$ 
\psp

Define
$$[u,v]=u\cdot v-(-1)^{i_1i_2}v\cdot u\qquad\for\;\;u\in\bar{\cal G}_i,\;v\in\bar{\cal G}_{i_2}.\eqno(4.62)$$
Pick any
$$L_0\in {\cal A}^-_0+{\cal A}^{(\iota+1)}_0\eqno(4.63)$$
(cf. (3.51) and (3.52)) and set
$$L=L_0+\sum_{i\in{Z}_2,\;j\in I_i}x_{i,j}\vs_{i,j}.\eqno(4.64)$$
Take any central element 
$$\kappa\in {\cal A}^-_0+{\cal A}^{(\iota+1)}_0\eqno(4.65)$$
of ${\cal A}$.  We define
$$\{f,g\}_1=\la (L,[\vs_{_{(f)}},(\kappa \vs_{_{(g)}})_+]-[(\kappa \vs_{_{(f)}})_-,\vs_{_{(g)}}]\ra\qquad\for\;\;f,g\in{\cal P}.\eqno(4.66)$$
Moreover, we define
$$\{f,g\}_2=\la (L\vs_{_{(f)}})_-L-L(\vs_{_{(f)}}L)_-,\vs_{_{(g)}}\ra\qquad\for\;\;f,g\in{\cal P}.\eqno(4.67)$$

The following is the  main theorem in this section.
\psp

{\bf Theorem 4.3}. {\it The brackets} $\{\cdot,\cdot\}_1$ {\it in (4.66) and} $\{\cdot,\cdot\}_2$ {\it in (4.67) forms a Poisson superpair on the algebra} ${\cal P}$.
\psp

{\it Proof}.  Let $\es\in\Bbb{F}$ be any fixed constant. Set
$$\hat{L}=L+\es\kappa.\eqno(4.68)$$
Then  we have
$$\ptl_u(\hat{L})=u\qquad\for\;\;u\in {\cal G}\eqno(4.69)$$
by (4.32) and (4.34). We define a map $H: \Omega\rta{\cal G}$ by
$$H=(\hat{L}u)_-\hat{L}-\hat{L}(u\hat{L})_-\qquad\for\;\;u\in{\cal G},\eqno(4.70)$$
For $\vs_{\eta}\in \Omega$, we have
$$H(\vs_{\eta})=\sum_{i\in\Bbb{Z}_2,\;j\in I_i}\eta_{i,j}((\hat{L}\vs_{i,j})_-\hat{L}-\hat{L}(\vs_{i,j}\hat{L})_-).\eqno(4.71)$$
Hence $H$ is of the form (4.45) by (4.63), (4.64) and the fact that $\kappa\in {\cal A}_0$. Moreover, it is straightforward to verify that
$$H(u)=(Lu)_-L-L(uL)_++\es(\kappa[L,u]_-+[(\kappa u)_-,L])\qquad\for\;\;u\in{\cal G}.\eqno(4.72)$$
Thus by (2.32)-(2.34), (4.49), (4.58), (4.66) and (4.67), we have
$$\{\cdot,\cdot\}_H=\es\{\cdot,\cdot\}_1+\{\cdot,\cdot\}_2.\eqno(4.73)$$
 Therefore, we only need to prove that $\{\cdot,\cdot\}_H$ is a Poisson superbracket on ${\cal P}$. Lemma 4.1 tells us that it is enough to prove (4.50) and (4.51).

For $u\in\Omega_{i_1}$ and $v\in\Omega_{i_2}$, we have
\begin{eqnarray*}\hspace{1cm}\la H(v),u\ra&=&\la (\hat{L}v)_-\hat{L}-\hat{L}(v\hat{L})_-,u\ra\\ &=&
\la [\hat{L}v-(\hat{L}v)_+]\hat{L}-\hat{L}[v\hat{L}-(v\hat{L})_+],u\ra\\ &=& -\la  (\hat{L}v)_+\hat{L}
-\hat{L}(v\hat{L})_+,u\ra\\ &=&
 -\la  (\hat{L}v)_+\hat{L},u\ra+\la \hat{L}(v\hat{L})_+,u\ra\\ &=&  -\la  (\hat{L}v)_+,\hat{L}u\ra+\la (v\hat{L})_+,u\hat{L}\ra\\&=& -\la \hat{L}v,(\hat{L}u)_-\ra+\la v\hat{L},(u\hat{L})_-\ra\\ &=& -(-1)^{i_1i_2}\la (\hat{L}u)_-,\hat{L}v\ra+\la v,\hat{L}(u\hat{L})_-\ra\\ &= & -(-1)^{i_1i_2}(\la (\hat{L}u)_-\hat{L},v\ra-\la \hat{L}(u\hat{L})_-,v\ra)\\ &=& -(-1)^{i_1i_2}\la (\hat{L}u)_-\hat{L}- \hat{L}(u\hat{L})_-,v\ra\\&=&-(-1)^{i_1i_2}\la H(u),v\ra\hspace{7.9cm}(4.74)\end{eqnarray*}
by (2.33), (2.34), (4.58) and the supersymmetry and associativity of $\la\cdot,\cdot\ra$ on $\bar{\cal G}$. Thus (4.50) holds.

Let  $u\in\Omega_{i_1},\;v\in\Omega_{i_2}$ and $w\in\Omega_{i_3}$. We have
\begin{eqnarray*}& &\la\ptl_{H(u)}(H)(v),w\ra\\ &=&\la (H(u)v)_-\hat{L}+(-1)^{i_1i_2}(\hat{L}v)_-H(u)-H(u)(v\hat{L})_--(-1)^{i_1i_2}\hat{L}(vH(u))_-,w\ra\\ &=& \la [((\hat{L}u)_-\hat{L}-\hat{L}(u\hat{L})_-)v]_-\hat{L}+(-1)^{i_1i_2}(\hat{L}v)_-((\hat{L}u)_-\hat{L}-\hat{L}(u\hat{L})_-)\\ && -((\hat{L}u)_-\hat{L}-\hat{L}(u\hat{L})_-)(v\hat{L})_--(-1)^{i_1i_2}\hat{L}[v((\hat{L}u)_-\hat{L}-\hat{L}(u\hat{L})_-)]_-,w\ra\\ &=& \la ((\hat{L}u)_-\hat{L}v)_-\hat{L},w\ra-\la (\hat{L}(u\hat{L})_-v)_-\hat{L},w\ra+(-1)^{i_1i_2}\la(\hat{L}v)_-(\hat{L}u)_-\hat{L},w\ra\\ & &-(-1)^{i_1i_2}\la(\hat{L}v)_-\hat{L}(u\hat{L})_-,w\ra
 -\la(\hat{L}u)_-\hat{L}(v\hat{L})_-,w\ra+\la\hat{L}(u\hat{L})_-(v\hat{L})_-,w\ra\\ & &-(-1)^{i_1i_2}\la \hat{L}(v(\hat{L}u)_-\hat{L})_-,w\ra+(-1)^{i_1i_2}\la\hat{L}(v\hat{L}(u\hat{L})_-)_-,w\ra\\ &=&\la ((\hat{L}u)_-\hat{L}v)_-,\hat{L}w\ra-\la (\hat{L}(u\hat{L})_-v)_-,\hat{L}w\ra+(-1)^{i_1i_2}\la(\hat{L}v)_-(\hat{L}u)_-,\hat{L}w\ra\\ & &-(-1)^{i_1i_2}\la(\hat{L}v)_-\hat{L}(u\hat{L})_-,w\ra-\la(\hat{L}u)_-\hat{L}(v\hat{L})_-,w\ra+\la(u\hat{L})_-(v\hat{L})_-,w\hat{L}\ra\\& &-(-1)^{i_1i_2}\la (v(\hat{L}u)_-\hat{L})_-,w\hat{L}\ra+(-1)^{i_1i_2}\la(v\hat{L}(u\hat{L})_-)_-,w\hat{L}\ra \\ &=&\la (\hat{L}u)_-\hat{L}v,(\hat{L}w)_+\ra-\la \hat{L}(u\hat{L})_-v,(\hat{L}w)_+\ra+(-1)^{i_1i_2}\la(\hat{L}v)_-(\hat{L}u)_-,\hat{L}w\ra\\& &-(-1)^{i_1i_2}\la(\hat{L}v)_-\hat{L}(u\hat{L})_-,w\ra-\la(\hat{L}u)_-\hat{L}(v\hat{L})_-,w\ra+\la(u\hat{L})_-(v\hat{L})_-,w\hat{L}\ra\\ & &-(-1)^{i_1i_2}\la v(\hat{L}u)_-\hat{L},(w\hat{L})_+\ra+(-1)^{i_1i_2}\la v\hat{L}(u\hat{L})_-,(w\hat{L})_+\ra\\ &=&(-1)^{i_1(i_2+i_3)}\la \hat{L}v, (\hat{L}w)_+ (\hat{L}u)_-\ra+(-1)^{i_1i_2}\la v\hat{L}, (u\hat{L})_-(w\hat{L})_+\ra\\ & &-[\la \hat{L}(u\hat{L})_-v,(\hat{L}w)_+\ra+\la(\hat{L}u)_-\hat{L}(v\hat{L})_-,w\ra-\la(u\hat{L})_-(v\hat{L})_-,w\hat{L}\ra]\\ & &-(-1)^{i_1i_2}[\la(\hat{L}v)_-\hat{L}(u\hat{L})_-,w\ra+\la v(\hat{L}u)_-\hat{L},(w\hat{L})_+\ra-\la(\hat{L}v)_-(\hat{L}u)_-,\hat{L}w\ra]\hspace{1.1cm}(4.75)\end{eqnarray*}
by (2.33), (2.34), (4.58) and the supersymmetry and associativity of $\la\cdot,\cdot\ra$ on $\bar{\cal G}$.
Moreover, 
\begin{eqnarray*}& &\la \hat{L}(u\hat{L})_-v,(\hat{L}w)_+\ra+\la(\hat{L}u)_-\hat{L}(v\hat{L})_-,w\ra-\la(u\hat{L})_-(v\hat{L})_-,w\hat{L}\ra+c.p.\\&=&
\la \hat{L}(u\hat{L})_-v,(\hat{L}w)_+\ra+(-1)^{(i_1+i_2)i_3}\la(\hat{L}w)_-\hat{L}(u\hat{L})_-,v\ra-\la(u\hat{L})_-(v\hat{L})_-,w\hat{L}\ra+c.p.\\&=&
\la \hat{L}(u\hat{L})_-v,(\hat{L}w)_+\ra+\la\hat{L}(u\hat{L})_-v,(\hat{L}w)_-\ra-\la(u\hat{L})_-(v\hat{L})_-,w\hat{L}\ra+c.p.\\&=& 
\la \hat{L}(u\hat{L})_-v,(\hat{L}w)_++(\hat{L}w)_-\ra-\la(u\hat{L})_-(v\hat{L})_-,w\hat{L}\ra+c.p.\\ &=& \la \hat{L}(u\hat{L})_-v,\hat{L}w\ra-\la(u\hat{L})_-(v\hat{L})_-,w\hat{L}\ra+c.p.\\ &=& \la (u\hat{L})_-v\hat{L},w\hat{L}\ra-\la(u\hat{L})_-(v\hat{L})_-,w\hat{L}\ra+c.p.
\\ &=&\la (u\hat{L})_-(v\hat{L})_+,w \hat{L}\ra+c.p.\\ &=& \la u\hat{L}v\hat{L},w \hat{L}\ra\hspace{11.8cm}(4.76)\end{eqnarray*}
and
\begin{eqnarray*}& &\la(\hat{L}v)_-\hat{L}(u\hat{L})_-,w\ra+\la v(\hat{L}u)_-\hat{L},(w\hat{L})_+\ra-\la(\hat{L}v)_-(\hat{L}u)_-,\hat{L}w\ra+c.p.\\&=&(-1)^{i_2(i_1+i_3)}\la(\hat{L}u)_-\hat{L}(w\hat{L})_-,v\ra+\la v(\hat{L}u)_-\hat{L},(w\hat{L})_+\ra-\la(\hat{L}v)_-(\hat{L}u)_-,\hat{L}w\ra+c.p.\\&=&\la v(\hat{L}u)_-\hat{L},(w\hat{L})_-\ra+\la v(\hat{L}u)_-\hat{L},(w\hat{L})_+\ra-\la(\hat{L}v)_-(\hat{L}u)_-,\hat{L}w\ra+c.p.\\&=&\la v(\hat{L}u)_-\hat{L},(w\hat{L})_-+(w\hat{L})_+\ra-\la(\hat{L}v)_-(\hat{L}u)_-,\hat{L}w\ra+c.p.\\&=&\la v(\hat{L}u)_-\hat{L},w\hat{L}\ra-\la(\hat{L}v)_-(\hat{L}u)_-,\hat{L}w\ra+c.p.\\&=&\la \hat{L}v(\hat{L}u)_-,\hat{L}w\ra-\la(\hat{L}v)_-(\hat{L}u)_-,\hat{L}w\ra+c.p.\\ &=&\la(\hat{L}v)_+(\hat{L}u)_-,\hat{L}w\ra+c.p.\\ &=& \la\hat{L}v\hat{L}u,\hat{L}w\ra\hspace{11.9cm}(4.77)\end{eqnarray*}
by (4.55)-(4.58) and the supersymmetry and associativity of $\la\cdot,\cdot\ra$ on $\bar{\cal G}$.
Thus
\begin{eqnarray*}& &\la\ptl_{H(u)}(H)(v),w\ra+c.p.\\ &=&[(-1)^{i_1(i_2+i_3)}\la \hat{L}v, (\hat{L}w)_+ (\hat{L}u)_-\ra+c.p.]+[(-1)^{i_1i_2}\la v\hat{L}, (u\hat{L})_-(w\hat{L})_+\ra+c.p.]\\ & &-[\la \hat{L}(u\hat{L})_-v,(\hat{L}w)_+\ra+\la(\hat{L}u)_-\hat{L}(v\hat{L})_-,w\ra-\la(u\hat{L})_-(v\hat{L})_-,w\hat{L}\ra+c.p.]\\& &-(-1)^{i_1i_2}[\la(\hat{L}v)_-\hat{L}(u\hat{L})_-,w\ra+\la v(\hat{L}u)_-\hat{L},(w\hat{L})_+\ra-\la(\hat{L}v)_-(\hat{L}u)_-,\hat{L}w\ra+c.p.]\\ &=&\la \hat{L}u,\hat{L}v\hat{L}w\ra+(-1)^{i_1i_2}\la v\hat{L}, u\hat{L}w\hat{L}\ra-\la u\hat{L}v\hat{L},w \hat{L}\ra-(-1)^{i_1i_2}\la\hat{L}v\hat{L}u,\hat{L}w\ra=0\hspace{0.5cm}(4.78)\end{eqnarray*}
by (4.56)-(4.58), the supersymmetry and associativity of $\la\cdot,\cdot\ra$ on $\bar{\cal G}$ and (4.75)-(4.77). Therefore, (4.51) holds.$\qquad\Box$
\psp

\section{Hamiltonian Superpairs in Variational Calculus}

In this section, we shall construct certain Hamiltonian superpairs in the formal variational calculus over a finite-dimensional $\Bbb{Z}_2$-graded associative algebras with a supersymmetric nondegenerate associative bilinear form. We assume that $\Bbb{F}$ is the field of real numbers or the field of complex numbers. 

Let ${\cal A}$ be a finite-dimensional $\Bbb{Z}_2$-graded associative algebra with a supersymmetric nondegenerate associative bilinear form $\la\cdot,\cdot\ra$. Denote
$$\dim {\cal A}_i=k_i\qquad\for\;\;i\in\Bbb{Z}_2.\eqno(5.1)$$
Fix a nonnegative intger $\iota$. Set
$$I_i=\ol{1,k_i}\times (-\Bbb{N}\bigcup\ol{1,\iota})\qquad\for\;\;i\in\Bbb{Z}_2.\eqno(5.2)$$
For $i\in\Bbb{Z}_2$, let $\{u_{i,j}\mid j\in I_i\}$
be the $C^{\infty}$-functions in real variable $x$,  taking values in ${\cal E}_i$ (cf. (4.2)). Set
$$\ptl={d\over dx},\;\;u_{i,j}^{(m)}={d^mu_{i,j}\over dx^m}\qquad\for\;\;m\in\Bbb{N},\;i\in\Bbb{Z}_2,\;j\in I_i.\eqno(5.3)$$
Denote by ${\cal P}$ the algebra supersymmetric polynomials in $\{u^{(m)}_{i,j}\mid i\in\Bbb{Z}_2,\; j\in I_i\}$. Then ${\cal P}$ has the $\Bbb{Z}_2$-grading
$${\cal P}_i=\mbox{Span}\:\{u^{(m_1)}_{i_1,j_1}\cdots u^{(m_\ell)}_{i_\ell,j_\ell}\mid \ell,m_p\in\Bbb{N},\;i_p\in\Bbb{Z}_2,\;j_p\in I_{i_p},\;\sum_{r=1}^\ell i_r\equiv i\}\eqno(5.4)$$
for $i\in\Bbb{Z}_2$, and (4.14) holds.

Take a basis $\{\vs_{i,l}\mid l\in\ol{1,k_i}\}$ of ${\cal A}_i$ for $i\in\Bbb{Z}_2$. Let $\hat{\cal A}$ be the free ${\cal P}$-module generated by $\{\vs_{i,l}\mid l\in\ol{1,k_i}\}$ with the $\Bbb{Z}_2$-grading
$$\hat{\cal A}_i=\sum_{j\in\ol{1,k_0}}{\cal P}_i\vs_{0,j}+\sum_{j\in\ol{1,k_1}}{\cal P}_{1+i}\vs_{1,j}\qquad\for\;\;i\in\Bbb{Z}_2.\eqno(5.5)$$
Moreover, we define the multiplication on $\hat{\cal A}$ by
$$(f\vs_{i_1,j_1})(g\vs_{i_2,j_2})=(-1)^{i_1p}fg\vs_{i_1,j_1}\vs_{i_2,j_2}\qquad\for\;\;i_q\in\Bbb{Z}_2,\;j_q\in\ol{1,k_{i_q}},\;g\in{\cal P}_p,\eqno(5.6)$$
and the bilinear form on $\hat{\cal A}$ by
$$\la f\vs_{i_1,j_1},g\vs_{i_2,j_2}\ra=(-1)^{i_1p}fg\la\vs_{i_1,j_1},\vs_{i_2,j_2}\ra\qquad\for\;\;i_q\in\Bbb{Z}_2,\;j_q\in\ol{1,k_{i_q}},\;g\in{\cal P}_p.\eqno(5.7)$$
Then $\hat{\cal A}$ forms a  $\Bbb{Z}_2$-graded associative algebra with a supersymmetric nondegenerate associative bilinear form $\la\cdot,\cdot\ra$. The algebra $\hat{\cal A}$ is an extension of ${\cal A}$ with extended supersymmetric bilinear form $\la\cdot,\cdot\ra$.

We view
$$\ptl=\sum_{i\in\Bbb{Z}_2,\;j\in\ol{1,k_i},\;m\in\Bbb{N}}u_{i,l}^{(m+1)}\ptl_{u_{i,j}^{(m)}}\eqno(5.8)$$
as a derivation of ${\cal P}$. Moreover, we extend $\ptl$ to a derivation of $\hat{\cal A}$ by
$$\ptl(f\vs_{i,j})=\ptl(f)\vs_{i,j}\qquad\for\;\;i\in\Bbb{Z}_2,\;j\in\ol{1,k_i}.\eqno(5.9)$$
Furthermore, we define the algebra of pseudo-differential operators
$${\cal D}=\{\sum_{l=-\infty}^n\phi_l\ptl^l\mid n\in\Bbb{Z},\;\phi_l\in\hat{\cal A}\}\eqno(5.10)$$
with multiplication determined by
$$(\phi\ptl^m)(\psi\ptl^n)=\sum_{p=0}^{\infty}\left(\!\!\begin{array}{c}m \\ p\end{array}\!\!\right)\phi\psi^{(p)}\ptl^{m+n-p}\qquad\for\;\;\phi,\psi\in\hat{\cal A},\;\;m,n\in\Bbb{Z},\eqno(5.11)$$
where 
$$\psi^{(p)}=\ptl^p(\psi).\eqno(5.12)$$
 Then ${\cal D}$ forms a $\Bbb{Z}_2$-graded associative algebra with the grading
$${\cal D}_i=\{\sum_{l=-\infty}^n\phi_l\ptl^l\mid n\in\Bbb{Z},\;\phi_l\in\hat{\cal A}_i\}\qquad\for\;\;i\in\Bbb{Z}_2.\eqno(5.13)$$

Define the space
$${\cal G}=\{\sum_{l=-\infty}^\iota\phi_l\ptl^l\mid\phi_l\in\hat{\cal A}\}.\eqno(5.14)$$
Then ${\cal G}$ is a $\Bbb{Z}_2$-graded subspace of ${\cal D}$ with the grading
$${\cal G}_0={\cal G}\bigcap {\cal D}_0,\;\;{\cal G}_1={\cal G}\bigcap {\cal D}_1.\eqno(5.15)$$
For
$$v=\sum_{i\in\Bbb{Z}_2,(j_1,j_2)\in I_i}f_{i,(j_1,j_2)}\vs_{i,j_1}\ptl^{j_2}\in {\cal G}\;\;\mbox{with}\;\;f_{i,(j_1,j_2)}\in {\cal P},\eqno(5.16)$$
we define the derivation of ${\cal P}$:
$$\ptl_v=\sum_{i\in\Bbb{Z}_2,\;j\in I_i,\;m\in\Bbb{N}}f^{(m)}_{i,j}\ptl_{u^{(m)}_{i,j}}\eqno(5.17)$$
(cf. (5.2)). It can be verified that as derivations on ${\cal P}$,
$$[\ptl,\ptl_v]=0\qquad\for\;\;v\in{\cal G}.\eqno(5.18)$$
Moreover, we set
$$\ptl_v(w)=\sum_{i\in\Bbb{Z}_2,(j_1,j_2)\in I_i}\ptl_v(g_{i,(j_1,j_2)})\vs_{i,j_1}\ptl^{j_2}\eqno(5.19)$$
for $v,w=\sum_{i\in\Bbb{Z}_2,(j_1,j_2)\in I_i}g_{i,(j_1,j_2)}\vs_{i,j_1}\ptl^{j_2}\in{\cal G}$. Furthermore, we can define a Lie superbracket on ${\cal G}$ by
$$[v,w]_0=\ptl_v(w)-(-1)^{i_1i_2}\ptl_w(v)\qquad\for\;\;v\in{\cal G}_{i_1},\;w\in{\cal G}_2.\eqno(5.20)$$

Set
$$\td{\cal P}={\cal P}/\ptl({\cal P})\eqno(5.21)$$
and denote
$$\td{f}=f+\ptl({\cal P})\qquad\for\;\;f\in{\cal P}.\eqno(5.22)$$
We define an action of ${\cal G}$ on $\td{\cal P}$ by
$$v(\td{f})=(\ptl_v(f))^{\sim}\qquad\for\;\;f\in{\cal P}.\eqno(5.23)$$
Then $\td{\cal P}$ becomes a ${\cal G}$-module. Moreover, we define a bilinear map $\la\cdot,\cdot\ra: {\cal D}\times{\cal D}\rta \td{\cal P}$ by
$$\la v,w\ra=\sum_{j\in\Bbb{Z}}(\la \phi_j,\psi_{-j-1}\ra)^{\sim}\qquad\for\;\;v=\sum_{m\in\Bbb{Z}}\phi_m\ptl^m,\;w=\sum_{m\in\Bbb{Z}}\psi_m\ptl^m\in {\cal D}\eqno(5.24)$$
(cf. (5.7)), where the sum is finite by (5.10). Then $\la\cdot,\cdot\ra$ forms a supersymmetric associative bilinear map (cf. (2.32)) by Lemma 2.6 in [DS].

For $v=\sum_{j\in\Bbb{Z}}\phi_j\ptl^j\in{\cal D}$, we define
and
$$v_+=\sum_{j=0}^{\infty}\phi_j\ptl^j,\qquad v_-=\sum_{j=1}^{\infty}\phi_{-j}\ptl^{-j}.\eqno(5.25)$$
Set
$${\cal D}_{\pm}=\{v_{\pm}\mid v\in{\cal D}\}.\eqno(5.26)$$
Then ${\cal D}^{\pm}$ are isotropic $\Bbb{Z}_2$-graded subalgebra of ${\cal D}$
with respect to the bilinear map $\la\cdot,\cdot\ra$. 
Define 
$$\Omega=\sum_{m=-\iota-1}^\infty\ptl^m\hat{\cal A}.\eqno(5.27)$$
Then $\Omega $ is a $\Bbb{Z}_2$-graded subspace of ${\cal D}$ by (5.11) with the grading
$$\Omega_0=\Omega\bigcap {\cal D}_0,\;\;\Omega_1=\Omega\bigcap {\cal D}_1.\eqno(5.28)$$
Identify $\Omega$ with a subspace of one-forms (taking values in $\td{\cal P}$) by
$$\varpi(v)=\la v,\varpi\ra\qquad\for\;\;v\in{\cal G},\;\varpi\in\Omega.\eqno(5.29)$$

We define {\it variational operators} on ${\cal P}$ by
$$\dlt_{(i,j)}=\sum_{m=0}^{\infty}(-\ptl)^m\circ\ptl_{u_{i,j}^{(m)}}\qquad\for\;\;i\in\Bbb{Z}_2,\;j\in I_i,\eqno(5.30)$$
where $\circ$ is the composition of operators on ${\cal P}$. It can be verified that
$$\dlt_{(i,j)}(\ptl({\cal P}))=\{0\}\qquad\for\;\;i\in\Bbb{Z}_2,\;j\in I_i.\eqno(5.31)$$
Moreover, we define a linear map $\chi: {\cal P}\rta \Omega$ by
$$\chi_{_f}=\sum_{i\in\Bbb{Z},\;j=(j_1,j_2)\in I_i}(-1)^{i(1+p)}\ptl^{-j_2-1}\dlt_{(i,j)}(f)\vs_{i,j_1}\qquad\for\;\;f\in{\cal P}_p.\eqno(5.32)$$
Then we can verify
$$df(v)=\la v,\chi_{_f}\ra\qquad\for\;\;f\in{\cal P},\;v\in{\cal G}.\eqno(5.33)$$
Hence
$$df=\chi_{_f}\in\Omega\qquad\for\;\;f\in{\cal P}.\eqno(5.34)$$

The Lie superbracket on ${\cal D}$ is defined by
$$[v_1,v_2]=v_1v_2-(-1)^{i_1i_2}v_2v_1\qquad\for\;\;v_1\in{\cal P}_{i_1},\;v_2\in{\cal P}_{i_2}.\eqno(5.35)$$
Take 
$$L_0=\sum_{m=-\infty}^{\iota+1}\sgm_m\ptl^m\;\;\mbox{with}\;\;\sgm_m\in{\cal A}_0\eqno(5.36)$$
and 
$$\kappa\in{\cal A}_0\bigcap(\mbox{Center}\:{\cal A}).\eqno(5.37)$$
Set
$$L=L_0+\sum_{i\in\Bbb{Z}_2,\;j=(j_1,j_2)\in I_i}u_{i,j}\vs_{i,j_1}\ptl^{j_2}\eqno(5.38)$$
(cf. (5.2)). We define the linear maps $H_1, H_2: \Omega\rta{\cal G}$ by
$$ H_1(\varpi)=\kappa[L,\varpi]_-+[(\kappa\varpi)_-,L],\;\;H_2(\varpi)=(L\varpi)_-L-L(\varpi L)_-\qquad\for\;\; \varpi\in\Omega.\eqno(5.39)$$

By the similar arguments as those in the proof of Theorem 4.3, we have:
\psp

{\bf Theorem 5.1}. {\it The pair} $(H_1,H_2)$ {\it forms a Hamiltonian superpair}.
\psp

{\bf Remark 5.2}. Suppose that ${\cal A}$ contains an identity element. Take the special case $L_0=\ptl^{\iota+1}$. Then
$$L=\ptl^{\iota+1}+\sum_{i\in\Bbb{Z}_2,\;j=(j_1,j_2)\in I_i}u_{i,j}\vs_{i,j_1}\ptl^{j_2}.\eqno(5.40)$$
By an algebraic manipulation, we can find
$$L^{1/(\iota+1)}=\ptl+\sum_{m=0}^{\infty}f_m\ptl^{-m}\;\;\;\mbox{with}\;\;f_m\in\hat{\cal A}\eqno(5.41)$$
such that
$$(L^{1/(\iota+1)})^{\iota+1}=L.\eqno(5.42)$$

Since ${\cal A}$ has an identity element $1_{\cal A}$, the map
$$\tr\:(a)=\la 1_{\cal A},a\ra \qquad\for\;\;a\in {\cal A}\eqno(5.43)$$
is a supersymmetric trace map, that is,
$$\tr\:(ab)=(-1)^{i_1i_2}\tr\:(ba)\qquad\for\;\;a\in {\cal A}_{i_1},\;b\in{\cal A}_{i_2}.\eqno(5.44)$$
We extent the trace map to ${\cal D}$ by
$$\tr\:(\sum_{i\in\Bbb{Z}_i,\;l\in\ol{1,k_i},\;m\in\Bbb{Z}}f_{i,l,m}\vs_{i,l}\ptl^m)=\sum_{i\in\Bbb{Z}_i,\;l\in\ol{1,k_i}}f_{i,l,-1}\tr\:(\vs_{i,l}).\eqno(5.45)$$
Then the above  map is supersymmetric trace map of the associative algebra ${\cal D}$ by Lemma 2.6 in [DS]. 

Assume that $\{u_{i,j}\mid i\in\Bbb{Z}_2,\;j\in I_i\}$ are also  $C^1$-functions of another variable $t$ and periodic in $x$. Take a positive integer $m$. Then the system
$${dL\over dt}=[L,(L^{m/(\iota+1)})_+]\eqno(5.46)$$
has infinitely many conservation laws:
$$\tr\:(L^{n/(\iota+1)})\qquad\for\;\;n\in\Bbb{N}+1.\eqno(5.47)$$

Suppose that $\{u_{i,j}\mid i\in\Bbb{Z}_2,\;j\in I_i\}$ are  $C^1$-functions of variable $\{t_1,t_2,t_3,....\}$ and periodic in $x$. Set
$$B_m=(L^{m/(\iota+1)})_+\qquad\for\;\;m\in\Bbb{N}+1.\eqno(5.48)$$
Assume that
$${dL\over dt_m}=[L,B_m]\qquad\for\;\;m\in\Bbb{N}+1.\eqno(5.49)$$
Then
$${dB_m\over dt_n}-{dB_n\over dt_m}=[B_m,B_n]\qquad\for\;\;m,n\in\Bbb{N}+1.\eqno(5.50)$$
The above equations are called the {\it equations of zero curvature}. We refer
to [DS], [D2], [Do] for more details.
\vspace{0.7cm}

\noindent{\Large \bf References}

\hspace{0.3cm}
\begin{description}

\item[{[A]}] M. Adler, On trace functional for formal pseudo-differential operators and the symplectic
structure of Korteweg-de Vries type equations, {\it Invent. Math.} {\bf 50} (1979), 219-248.

\item[{[BDH]}] J. C. Brunelli, A. Das and W.-J. Huang, Gelfand-Dikii brackets for nonstandard Lax equations, {\it Modern Phys. Lett.} {\bf A 9} (1994), 2147-2155.

\item[{[DH]}]A. Das and W.-J. Huang, The Hamiltonian structures associated with a generalized Lax operator, {\it J. Math. Phys.} {\bf 33} (1992), 2487-2497.

\item[{[DHP]}] A. Das, W.-J. Huang and S. Panda, The Hamiltonian structures of the KP hierarchy, {\it Phys. Lett.} {\bf B 271} (1991), 109-115. 

\item[{[DP]}] A. Das and S. Panda, Gelfand-Dikii brackets for nonstandard supersymmetric system, {\it Modern Phys. Lett.} {\bf A 11} (1996), 723-730.

\item[{[D1]}] L. A. Dickey, On Hamiltonian and Lagrangian formalisms for the KP-hierarchy of integrable equations, {\it Ann. N. Y. Acad. Sci.} {\bf 491} (1987), 131-148.

\item[{[D2]}] L. A. Dickey, {\it Soliton Equations and Hamiltonian Systems}, Advanced Series in Mathematical Physics, Vol. {\bf 12}, World Scientific. Singapore/New Jersey/London/Hong Kong, 1991.

\item[{[Do]}] I. Ya.  Dorfman, {\it Dirac structures and integrability of nonlinear evolution equations}, John Wiley \& Sons, Chichester/New York/Brisbane/Toronto/Singapore, 1993.

\item[{[DF]}] I. Ya.  Dorfman and A. S. Fokas, Hamiltonian theory over noncommutative rings and integrability in multidimensions, {\it J. Math. Phys.} {\bf 33} (1992), 2504-2514.

\item[{[DS]}] V. G. Drinfel'd and V. V. Sokolov, Lie algebras and equations ofKorteweg-de Vries type,
{\it J. Sov. Math.} {\bf 30} (1985), 1975-2035.

\item[{[GDi1]}] I. M. Gel'fand and L. A. Dikii, Asymptotic behaviour of the resolvent of Sturm-Liouville equations and the algebra of the Korteweg-de Vries equations, {\it Russian Math. Surveys} {\bf 30:5} (1975), 77-113.

\item[{[GDi2]}] I. M. Gel'fand and L. A. Dikii, A Lie algebra structure in a formal variational Calculation, {\it Func. Anal. Appl.}  {\bf 10} (1976), 16-22.

\item[{[GDi3]}] I. M. Gel'fand and L. A. Dikii, A family of Hamiltonian structures connected with integrable nonlinear partial differential equations, {\it Preprint No.} {\bf 136} {\it IPM AN SSSR, Moscow,} 1978.

\item[{[GDo1]}] 
I. M. Gel'fand and I. Ya. Dorfman, Hamiltonian operators and algebraic structures related to them, {\it Funkts. Anal. Prilozhen}  {\bf 13} (1979), 13-30.

\item[{[GDo2]}] I. M. Gel'fand and I. Ya. Dorfman, Schouten bracket and Hamiltonian operators, {\it Funkts. Anal. Prilozhen}  {\bf 14} (1980), 71-74.

\item[{[HKW]}] P. Harpe, M. kervaire and C. Weber, On the Jones polynomial, {\it L'Enseig. Math.} {\bf 32} (1986), 271-235.

\item[{[MR]}] Yu. I. Manin and A. O. Radul, A supersymmetric extension of the Kadomtsev-Petviashvili
hierarchy, {\it Commun. Math. Phys.} {\bf 98} (1985), 65-77.

\item[{[R]}]  A. O. Radul, Two series of Hamiltonian structures for the hierarchy of Kadomtsev-Petviashvili equations, in ``Applied methods of nonlinear analysis and control,'' ed. Mironov, Moroz, and Tshernjatin, MGU, (1987) (Russian), pp. 149-157.

\item[{[W]}] Y. Watanabe, Hamiltonian structure of Sato's hierarchy of KP equations and a coadjoint orbits of a certain formal Lie group, {\it Lett. Math. Phys.} {\bf 7} (1983), 99-106.

\item[{[X1]}] X. Xu, Hamiltonian superoperators, {\it J. Phys A: Math.  Gen.} {\bf 28}  (1995), 1681-1698.

\item[{[X2]}] X. Xu, Variational calculus of supervariables and related algebraic structures, {\it J. Algebra} {\bf 223} (2000), 396-437.

\item[{[X3]}] X. Xu, Equivalence of conformal superalgebras to Hamiltonian superoperators, to appear in {\it Algebra Colloquium} {\bf 8} (2001), 63-92.

\end{description}
\end{document}